%% file: main.tex
\documentclass[a4paper]{article}

\usepackage{booktabs}
\usepackage[english]{babel}
\usepackage[utf8]{inputenc}
\usepackage{amssymb}
\usepackage[fleqn]{amsmath}
\usepackage{amsthm}
\usepackage{graphicx}
\usepackage{framed}
\usepackage[colorinlistoftodos]{todonotes}
\usepackage{subcaption}
\usepackage{multirow}
\usepackage{enumitem}
\usepackage{titling} 
\usepackage{titlesec}
\usepackage[top=1.25in, bottom=1.25in, left=1.25in, right=1.25in]{geometry}
\usepackage{breqn}
\usepackage[numbers]{natbib}
\usepackage{tikz}
\usepackage{stackengine}
\usetikzlibrary{arrows,%
	petri,%
	topaths}%
\usepackage{tkz-berge}
\usepackage{systeme}
\usepackage{algorithm}
\usepackage{algpseudocode}
\usepackage{tabstackengine}
\def\typesystem#1#2{\savestack{#1}{\setstackEOL{,}\setstackTAB{ }
  $\left\{\ensurestackMath{\tabbedCenterstack[r]{#2}}\right.$}}
\TABbinary
\include{commands}
\stackMath
\newcommand{\trace}[1]{\text{Tr }#1}

\newcommand{\re}[1]{\text{Re}\left(#1\right)}
\newcommand{\rank}[1]{\text{Rank}\left(#1\right)}

\theoremstyle{plain}
\newtheorem{theorem}{Theorem}[section]

\theoremstyle{plain}
\newtheorem{lemma}[theorem]{Lemma}

\theoremstyle{plain}
\newtheorem{proposition}[theorem]{Proposition}

\theoremstyle{plain}
\newtheorem{corollary}{Corollary}[theorem]

\theoremstyle{definition}
\newtheorem{definition}{Definition}[section]

\theoremstyle{plain}
\newtheorem{conjecture}[theorem]{Conjecture}

\ERPWtemplate{Wissing \& Van Dam (2019)}{Strongly HDS digraphs}{}

\title{
The negative tetrahedron and the first infinite family of connected digraphs that are strongly determined by the Hermitian spectrum
}
\author{
Pepijn Wissing\thanks{Corresponding author: p.wissing@tilburguniversity.edu},~ Edwin R. van Dam \\ \small{CentER and Department of Econometrics and Operations Research, Tilburg University}
}
\date{}
\begin{document}
\linedabstractkw{
Thus far, digraphs that are uniquely determined by their Hermitian spectra have proven elusive. 
Instead, researchers have turned to spectral determination of classes of switching equivalent digraphs, rather than individual digraphs. 
In the present paper, we consider the traditional notion: a digraph (or mixed graph) is said to be strongly determined by its Hermitian spectrum (abbreviated SHDS) if it is isomorphic to each digraph to which it is cospectral. Convincing numerical evidence to support the claim that this property is extremely rare is provided. 
Nonetheless, the first infinite family of connected digraphs that is SHDS is constructed. This family is obtained via the introduction of twin vertices into a  structure that is named \textit{negative tetrahedron}. This special digraph, that exhibits extreme spectral behavior, is contained in the surprisingly small collection of all digraphs with exactly one negative eigenvalue, which is determined as an intermediate result. 
}{
Hermitian adjacency matrix, Spectra of digraphs, Cospectral digraphs, Spectral characterization of digraphs
}
\maketitle
\section{Introduction}
The relation between the eigenvalues of a graph and its structural characteristics has been studied by many. Finding use in various fields such as combinatorics \cite{brouwer2011spectra}, optimization \cite{lovasz1979shannon, mohar1993eigenvalues} and computer science \cite{spielman2007spectral}, the literature on graph spectra is vast. Of particular interest to the authors is the question: can we determine a graph by its spectrum? This question has received considerable attention; the known results have been surveyed in \citep{vandam2003, vandam2009}. 

Nevertheless, the same question has yielded far fewer results in the directed graph (or \textit{mixed graph}, see Sec. 1.1) paradigm. This is in part due to the absence of a consensus as to which matrix best reflects the characteristics of a directed graph (henceforth abbreviated as \textit{digraph}) in its eigenvalues. 
A natural choice is the well-known \textit{adjacency matrix}. 
However, as this matrix is in general not symmetric, its eigenvalues are not necessarily real. 
Other alternatives include the \textit{skew-symmetric adjacency matrix} (see \citep{cavers2012skew}), which only works for digraphs that contain no bidirected edges, and the \textit{normalized Laplacian} \citep{chung2005laplacians}, which has a particularly impractical definition.

The introduction of the Hermitian adjacency matrix ("Hermitian" or "$H$", for short), independently by \citet{liu2015hermitian} and \citet{guo2016hermitian}, has provided us with an interesting new candidate. A number of important pieces of machinery from undirected graph theory have been shown to hold with respect to the Hermitian. These include, but are not limited to, eigenvalue interlacing and the quotient matrix \citep{haemers1995, guo2016hermitian} and the concept of computing the coefficients of the characteristic polynomial by considering elementary subdigraphs \citep{liu2015hermitian}. 

One of the main drawbacks of the Hermitian, from a spectral analysis point of view, is that a large number of digraphs with the same underlying graph have an identical $H$-spectrum. Suppose, for example, that a digraph has a cut-edge. Then this edge may be directed in either direction, or even be made bidirected, all the while leaving the corresponding $H$-spectrum unchanged. This illustration is a special case of an operation uncovered by \citet{guo2016hermitian}, called \textit{four-way switching}, which partitions the digraph in four (possibly empty) parts and performs a series of operations that either change the direction or nature (from \textit{arc} to \textit{digon} or vice versa) of the arcs between the parts or take the converse of the digraph as a whole, without affecting the $H$-spectrum. 

In recent days, the Hermitian spectrum has been the subject of several publications. 
In addition to co-authoring one of the two fundamental papers \citep{liu2015hermitian, guo2016hermitian}, mentioned above, \citet{mohar2016hermitian} has characterized all digraphs whose Hermitian adjacency matrix have rank 2, and shown that there are infinitely many digraphs that have cospectral mates, which are not members of the same switching equivalence class.
\citet{wang2017mixed} extends the research in \citep{mohar2016hermitian} to the digraphs of rank 3; their main result is that any pair of weakly connected, cospectral rank 3 digraphs is switching equivalent. 
Although using a different approach, \citet{tian2018nullity} obtain similar results as in \citep{wang2017mixed}.

Further recent research that is concerned with the Hermitian spectrum but less relevant to the current paper includes  \citet{greaves2019interlacing}, \citet{guo2017digraphs}, \citet{greaves2012cyclotomic}, \citet{hu2017spectral}, \citet{chen2018relationmatching} and \citet{chen2018relationrank}.\\

In his article \citep{mohar2016hermitian}, Mohar investigates which digraphs of rank 2 are determined by their $H$-spectrum. (Abbreviated HDS hereafter.) Due to the before-mentioned characteristics of the Hermitian, Mohar defines a digraph to be HDS if it is cospectral only to those digraphs that are obtained from the digraph by a switching operation, possibly followed by the reversal of all edges. 
This definition is, however, much weaker than that of the similarly named concept in undirected graph context; if a graph $G$ is said to be determined by its adjacency spectrum, then one is able to uniquely (up to isomorphism) construct said graph when one is given its spectrum. 

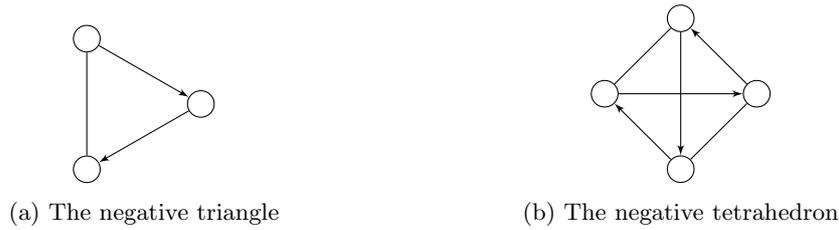
\begin{figure}[h!]
\begin{center}
\begin{subfigure}[b]{0.45\textwidth}\centering
\begin{tikzpicture} 
\tikzset{vertex/.style = {shape=circle,draw,minimum size=1em}}
\tikzset{edge/.style = {->,> = latex'}}
\tikzset{symedge/.style = {-,> = latex'}}
\tikzset{edgeBendLeft/.style = {->,> = latex',bend left}}
\tikzset{edgeBendRight/.style = {->,> = latex',bend right}}


\def\n{3}
\def\rad{1}
\foreach \a in {1,2,...,\n}{
\node[vertex] (\a) at (\a*360/\n: \rad cm) {};
}
\draw[symedge] (1) to node{} (2);
\draw[edge] (1) to node{} (3);
\draw[edge] (3) to node{} (2);


\end{tikzpicture}
\caption{The negative triangle}
\label{fig: example SHDS digraph a}

\end{subfigure}
\quad
\begin{subfigure}[b]{0.45\textwidth}\centering
\begin{tikzpicture} 
\tikzset{vertex/.style = {shape=circle,draw,minimum size=1em}}
\tikzset{edge/.style = {->,> = latex'}}
\tikzset{symedge/.style = {-,> = latex'}}
\tikzset{edgeBendLeft/.style = {->,> = latex',bend left}}
\tikzset{edgeBendRight/.style = {->,> = latex',bend right}}


\def\n{4}
\def\rad{1}
\foreach \a in {1,2,...,\n}{
\node[vertex] (\a) at (\a*360/\n: \rad cm) {};
}
\draw[symedge] (1) to node{} (2);
\draw[symedge] (3) to node{} (4);
\draw[edge] (1) to node{} (3);
\draw[edge] (3) to node{} (2);
\draw[edge] (2) to node{} (4);
\draw[edge] (4) to node{} (1);


\end{tikzpicture}
\caption{The negative tetrahedron}
\label{fig: example SHDS digraph b}

\end{subfigure}

\end{center}
\caption{\label{fig: example SHDS digraph}Two strongly HDS digraphs.}
\end{figure}

As such, the authors set out to classify digraphs which are \textit{strongly} determined by their $H$-spectrum; that is, digraphs whose spectra occur uniquely. 
Two prominent examples of such digraphs are shown in Figure \ref{fig: example SHDS digraph}.
These digraphs are extremely rare, as any such digraph must be self-converse.
We observe that the fraction of digraphs that satisfies this property rapidly goes to zero as the number of vertices grows. 
While no formal proof was found in the literature, numerical evaluation of the counting polynomials by \citet{harary1955number} and \citet{harary1966enumeration} provides convincing evidence to support this claim.

In the present paper, the authors were inspired by a result first encountered in \cite{guo2015simple} and later in \cite{guo2016hermitian}, which occurs here as Lemma \ref{lemma: contained triangle H3}. 
In particular, by said lemma, there is exactly one kind of induced subdigraph of order 3 that may contribute negatively to $\trace{H(D)^3}$, where $H(D)$ is the Hermitian of a digraph $D$.  
This order $3$ digraph is named \textit{negative triangle} and shown in Figure \ref{fig: example SHDS digraph a}. 
Furthermore, eigenvalue interlacing is used extensively.

The main results of this paper are as follows. We construct the first infinite family of connected, strongly HDS digraphs, in Theorem \ref{thm: main theorem}. This family is obtained by twin expansion (see Def. \ref{def: twin expansion}) of a key digraph, which is named the \textit{negative tetrahedron.} 
This peculiar digraph is a tetrahedron, whose four faces are all negative triangles, as is shown in Figure \ref{fig: example SHDS digraph b}. Moreover, it is the only  reduced\footnote{A digraph is said to be reduced if it contains no two vertices whose neighborhoods exactly coincide, and no isolated vertices. See Def. \ref{def: reduced}.} digraph that has rank 4 and exactly one negative eigenvalue.
Additionally, we determine all digraphs that have precisely one negative eigenvalue in Theorem \ref{thm: TR of a digraph with 1 neg ev}, and show that any pair of connected rank  $4$ members of this class is switching equivalent if they are cospectral in Proposition \ref{prop: connected cospectral rank 4 digraphs are sw eq}.

\section{Preliminaries}
Let us first thoroughly discuss the definitions and known results. Most of the concepts below are well-known and will not always be explicitly referenced.  For more detail, the reader is referred to e.g., \citet{bondy2008graph} or any other recent book on graph theory.
\subsection{Terminology and notation}
\label{sec: terminology and notation}
Throughout this work, we will adopt the following terminology and notation.
A \textit{directed graph} or \textit{digraph} of order $n$ is denoted $D=(V,E)$.
Here, we denote the vertex set of $D$ as $V$ (sometimes specified as $V(D)$) and edge set $E$ (or $E(D)$), where $E \subseteq V \times V$. 
An \textit{arc} or \textit{directed edge} $(u,v)\in E$ is an ordered pair of vertices; $u$ is called the \textit{initial vertex} and $v$ is called the \textit{terminal vertex}. 
If both $(u,v)\in E$ and $(v,u)\in E$, we say that the unordered pair $\{u,v\}$ is a \textit{digon} in $D$. 
A \textit{loop} is an arc of which the terminal vertex equals the initial vertex. 
Throughout, we will not allow digraphs that are considered in this work to contain loops.

Some of the cited literature considers \textit{mixed graphs}, which are defined to be an ordered triple $(V,E,A),$ where $V$ is the vertex set, $E$ is the undirected edge set and $A$ is the directed edge set. Since a single bidirected edge is, for our purposes throughout this paper, equivalent to two arcs whose directions are reversed, we consider the class of mixed graphs equivalent to the class of digraphs as described above. In the interest of clearness, we will only use the word "graph" when it concerns an \textit{undirected}  graph. 

The Hermitian adjacency matrix, introduced independently by \citet{guo2016hermitian} and \citet{liu2015hermitian}, is formally defined as follows. 
\begin{definition}\label{def: hermitian}\cite{guo2016hermitian, liu2015hermitian}
Let $D=(V,E)$ be a digraph of order $n$. Define the \textit{Hermitian} $H=H(D)$ as the $n\times n$ matrix with entries
\begin{equation}
\label{eq: def hermitian}
h_{uv} = \begin{cases}
1 & \text{ if } (u,v)\in E \text{ and } (v,u)\in E,\\
i & \text{ if } (u,v)\in E \text{ and } (v,u)\not\in E,\\
-i & \text{ if } (u,v)\not\in E \text{ and } (v,u)\in E,\\
0 & \text{ otherwise.}
\end{cases}
\end{equation}
\end{definition}
The Hermitian $H$ is diagonalizable \cite{brouwer2011spectra} and one may apply eigenvalue interlacing (see Lemma \ref{lemma: interlacing}) with respect to $H$, as shown in \cite{guo2015simple, guo2016hermitian}.
Throughout, we will use the Hermitian to characterize digraphs; the rank of $D$ and the spectrum of $D$ are said to be respectively the rank and the spectrum of the Hermitian $H(D)$ of $D$. 
Accordingly, the characteristic polynomial $\chi_D$ of $D$ is defined as $\chi_D(\mu) = \det\left(\mu I - H(D)\right)$. 
The spectrum of $D$ is denoted $\Sigma_D = \left\{\lambda_1^{[m_1]},\ldots,\lambda_k^{[m_k]}\right\}$, where $\lambda_1,\ldots,\lambda_k$ are the $k$ distinct eigenvalues, whose multiplicities are $m_1,\ldots,m_k$. 
As is convention, $\rho(D)$ denotes the spectral radius of $D$, i.e., its largest eigenvalue in absolute value. 

We may sometimes be interested in the \textit{underlying graph} of a given digraph $D$. 
Let $\Gamma(\cdot)$ be the operator that transforms a digraph into its underlying graph. 
That is, given $D$, the graph $\Gamma(D)$ is obtained by replacing every arc in $D$ with the corresponding digon. 
Two vertices are said to be adjacent in $D$ if they are adjacent in $\Gamma(D)$.
Moreover, note that we may determine the number of edges in the underlying graph of a given digraph from its spectrum by using that $\text{Tr}~H(D)^2 = 2\cdot |E(\Gamma(D))|$. 

The \textit{converse} of a digraph $D$ is denoted $D^\top$, and is obtained by reversing the direction of all arcs. Note that $H(D)^\top = H\left(D^\top\right)$, and thus any digraph is inherently cospectral to its converse. Furthermore, a digraph that is isomorphic to its converse is called \textit{self-converse}; this rare property is necessary for strong spectral determination. 

The last mentioned piece of notation is concerned with induced subdigraphs. If $D = (V,E)$ and let $W\subset V$, then we denote the induced subgraph that is obtained by removing any vertices in $V\setminus W$ and removing any edges that are incident to a vertex in $V\setminus W$ as $D[W]$.

As was briefly mentioned before, \citet{mohar2016hermitian} pioneered spectral characterization with respect to the Hermitian adjacency matrix. 
We include the details of Mohar's definitions below.
\begin{definition}\cite{guo2016hermitian}\label{def: 4 way switching}
A \textit{four-way switching} is the operation of changing a digraph $D$ into the digraph $D'$ by choosing an appropriate diagonal matrix $S$ with $S_{jj} \in \{\pm 1, \pm i\},~j = 1,\ldots,|V(D)|,$ and setting $H(D') = S^{-1}H(D)S.$ 
Informally, $S$ is appropriate when $H(D')$ is a Hermitian adjacency matrix. 
\end{definition}
\begin{definition}\label{def: switching equivalent}\cite{mohar2016hermitian}
Two digraphs $D$ and $D'$ are said to be \textit{switching equivalent} if $D'$ may be obtained from $D$ by a four-way switching, possibly followed by a reversal of all arcs.
\end{definition}

\begin{definition}\label{def: whds}
\citep{mohar2016hermitian} A digraph is said to be \textit{(weakly) determined by its Hermitian spectrum (WHDS)} if it is switching equivalent to every digraph to which it is $H$-cospectral.
\end{definition}
However, the following definition of spectral determination is in a sense more loyal to its undirected graph analogue. 
\begin{definition}\label{def: shds}
A digraph is said to be \textit{strongly determined by its Hermitian spectrum (SHDS)} if it is isomorphic to every digraph to which it is $H$-cospectral.
\end{definition}
The distinction between Definitions \ref{def: whds} and \ref{def: shds} summarizes what sets this work apart from previous articles regarding Hermitian spectral characterization. 
To distinguish between the two definitions, the authors have added the word "weakly" to the former and "strongly" to the latter. The terminology is justified by the observation that any SHDS digraph is implicitly WHDS. 

We end this section with a few words of warning, regarding a frequent mistake surrounding the term (W)HDS. Since neither four-way switching nor taking the converse changes a digraph's underlying graph, two digraphs of which just one is connected cannot be switching equivalent. 
Hence, if a connected digraph is cospectral to a digraph that contains isolated vertices, neither may be said to be (W)HDS. 
While \cite{wang2017mixed, tian2018nullity} show that any pair of cospectral, connected rank 3 digraphs is switching equivalent, cases in which such a connected digraph is cospectral to a disconnected digraph may still be found. 
Therefore, the phrasing of the final theorems of both \cite{wang2017mixed} and \cite{tian2018nullity}, in which it is claimed that any rank 3 digraph is (W)HDS, is flawed.

\subsection{Known results}
Here, we will list some of the results that are vital to the discussion in this paper. Likely the single most used tool throughout is known as eigenvalue interlacing (\citet{haemers1995}), which is a particularly powerful tool, adopted from graph theory. In \citep{godsil2001algebraic}, this concept has been shown to be valid with respect to Hermitian matrices. 
\begin{lemma}\label{lemma: interlacing} \citep{haemers1995, godsil2001algebraic}
Suppose $A$ is a Hermitian $n\times n$ matrix with eigenvalues $\lambda_1\geq \ldots\geq \lambda_n$. Then the eigenvalues $\mu_1\geq \ldots\geq \mu_m$ of a principal submatrix of size $m$ satisfy $\lambda_i \geq \mu_i \geq \lambda_{n-m+i}$ for $i\in[m]$.
\end{lemma}
In particular, Lemma \ref{lemma: interlacing} will be one of the tools to determine the collection of all digraphs with a single negative eigenvalue. To explicitly write the spectra of these families, the concepts known as \textit{quotient matrix} and \textit{equitable partition} are used. Both originate from \citep{haemers1995}, in graph context, and have been published in \citep{guo2016hermitian} for Hermitian context.

Let $D$ be a digraph and let $\mathcal{V} = \{V_1,\ldots,V_k\}$ be a partition of $V(D)$. One may order the vertices of $V$ such that $\mathcal{V}$ induces a partition of $H$ into block matrices as 
\[ H = \begin{bmatrix}
H_{11} & \cdots & H_{1k} \\ \vdots & \ddots & \vdots \\ H_{k1} & \cdots & H_{kk}
\end{bmatrix}.\]
The quotient matrix of a Hermitian $H$ with respect to $\mathcal{V}$ is the matrix $Q = [q_{ij}],$ $i,j \in [k]$, where $q_{ij}$ is the average row sum of block $H_{ij}$. The partition $\mathcal{V}$ is said to be \textit{equitable} if every block $H_{ij}$ has constant row sum. One then has the following result. 
\begin{lemma}\label{lemma: equitable partition quotient}\citep{guo2016hermitian}
Let $D$ be a digraph with Hermitian $H$, and let $\mathcal{V}$ be an equitable partition of its vertices. Moreover, let $Q$ be the quotient matrix of $H$ with respect to $\mathcal{V}.$ If $\lambda$ is an eigenvalue of $Q$ with multiplicity $\mu$, then it is an eigenvalue of $H$ with multiplicity at least $\mu$.
\end{lemma}

In \cite{mohar2016hermitian}, Mohar works with digraphs of rank 2, particularly those that are complete bipartite or complete tripartite. With regard to the former, the following unproven claim is made. 
\begin{conjecture}\citep{mohar2016hermitian}
There are only finitely many integers $m$ and $n$ for which the complete bipartite graph $K_{m,n}$ is WHDS. 
\end{conjecture}
Regarding complete tripartite digraphs, Mohar claims that there are many instances that are WHDS for number theoretic reasons. In particular, if we denote by $\overrightarrow{C}_3(a,b,c)$ the complete tripartite digraph with parts $A$, $B$, $C$, where $|A|=a$, $|B|=b$ and $|C| = c$, with all arcs from $A$ to $B$, $B$ to $C$, and $C$ to $A$, then the following claims hold.
\begin{proposition}\citep{mohar2016hermitian}\label{prop: mohar characterization}
 $\overrightarrow{C}_3(n,n,n)$, $\overrightarrow{C}_3(n,n,n+1)$ and $\overrightarrow{C}_3(n-1,n,n)$ are WHDS for every $n$. 
\end{proposition}
\begin{corollary}\cite{mohar2016hermitian}
Suppose that $a$ and $n>a>0$ are integers such that $a^2<2n.$ Then $\overrightarrow{C}_3(n-a,n,n+a)$ is WHDS if and only if $a$ is not divisible by a prime that is congruent to $1$ modulo $6$.
\end{corollary}

This line of research was extended to rank 3 independently by \citet{wang2017mixed} and \citet{tian2018nullity}. Of most relevance to this work is the following result.\footnote{Originally formulated as \textit{"All connected digraphs of order $n$ with rank $3$ are WHDS,"} which is not entirely accurate, as discussed at the end of Sec. \ref{sec: terminology and notation}. The authors took the liberty of slightly rephrasing to avoid confusion.} 
\begin{proposition}\label{prop: wang characterization}\cite{wang2017mixed, tian2018nullity}
Any two connected, cospectral, rank $3$ digraphs are switching equivalent.
\end{proposition}
However, \citet{wang2017mixed} also show that, if the assumption on connectedness is omitted, one may construct infinite families of rank 3 digraphs that are not WHDS.
\begin{proposition}\label{prop: wang characterizaton}\cite{wang2017mixed}
There are infinitely many digraphs with rank $3$ that are not WHDS. 
\end{proposition}

\subsection{Twins}
The first half of the discussion in this paper will concern 'small' digraphs with exactly one negative eigenvalue; the second half will extend this discussion to 'large' digraphs. That said, the discussed digraphs remain largely similar, from a structural point of view. Specifically, in order to increase the size of the considered digraphs without compromising the structural arguments made in the former part, \textit{twins} are introduced into the small digraphs. Since there have been several authors (e.g., \cite{axenovich2013regularity, butler2014using}) to have introduced a similarly named concept, each with subtle differences, we provide the formal definition as it is used throughout this paper.

\begin{definition}
Two vertices $u,v$ in $D$, whose Hermitian is $H$, are called \textit{twins} if $H_{ux} = H_{vx}$ for every $x\in V(D).$
\end{definition}
A simple but important observation is that $u$ and $v$ are implicitly not twins if $H_{uv}\not=0$, as loops are not allowed throughout. Moreover, if $u,v,w$ are vertices in $D$, $w$ is said to \textit{distinguish} $u$ and $v$ if $H_{uw} \not= H_{vw}$.  
Naturally, if such a vertex $w$ exists in $D$, then $u$ and $v$ are not twins in $D$, which justifies the terminology.

We will often want to consider the structure that is in a sense fundamental to a large digraph that contains a set of equivalent vertices. 
To this end, we define the \textit{twin reduction} operation, which reduces such a collection of twin vertices down to a single vertex; this may significantly reduce the order of a digraph, while retaining its general structure and rank.
Moreover, using said operation, we define a property that characterizes digraphs we consider to be 'small'. 
\begin{definition}
We define $TR(\cdot)$ to be the \textit{twin reduction} operator, which removes vertices and edges from a digraph in such a way that exactly one of every collection of twins is kept and no isolated vertices remain. 
\end{definition}
\begin{definition}\label{def: reduced}
A digraph $D$ is called \textit{reduced} if $TR(D) = D$.
\end{definition}

Naturally, we may also want to reverse this operation, to increase the size of the digraph without compromising the fundamental structure. The formal definition of the corresponding operation is given below. 
\begin{definition}
A vector $t = \begin{bmatrix}t_0 & t_1 & \cdots & t_n \end{bmatrix}\in\mathbb{N}_0^n$ is called an \textit{expansion vector} for a digraph $D$ if $n = |V(D)|$ and $t_1,\ldots,t_n \geq 1$. 
\end{definition}

\begin{definition}\label{def: twin expansion}
Let $D$ be a digraph with an ordered set $V$ of $n$ vertices, and let $t$ be an appropriate expansion vector.
The \textit{twin expansion of $D$ with respect to $t$} is denoted $D'=TE(D,t)$ and is obtained by replacing each vertex $u$ in $D$ by $t_u$ twins, and adding $t_0$ isolated vertices.
Formally, if $V=[n]$, let $V(D')=V_0 \cup V_1 \cup \cdots \cup V_n$, where $V_0,V_1,\dots,V_n$ are mutually disjoint sets, with $|V_u|=t_u$.
In $D'$, $V_0$ is a set of isolated vertices, and $H'_{u'v'} = H_{uv}$ for every $u'\in V_u$, $v'\in V_v$, $u,v\in V$, where $H$ and $H'$ are the Hermitians of $D$ and $D'$, respectively.
\end{definition}

Note that each entry of the expansion vector thus corresponds to a specific vertex in the digraph that is to be expanded. As a direct consequence, the vertex ordering does matter, in the sense that permuting the expansion vector does not change the resulting (expanded) digraph if the vertex order of the source digraph is permuted accordingly. Thus, we will fix the vertex orderings of the relevant digraphs to ensure that the above does not occur, when permutations of given expansion vectors are discussed. Throughout, the northmost vertex of any digraph defined through an illustration is assigned '$1$', after which vertices are subsequently labeled in counterclockwise order.

For the sake of clarity, we include the following example that shows the working of the twin expansion operator.

\example{ex1}{
Let $t = \begin{bmatrix} 2 & 3 & 2 & 1\end{bmatrix}.$ Then the vertices of the negative triangle $T^-$ and its twin expansion $D' = TE(T^-, t)$ may be labeled such that their Hermitians are
\[H(T^-) = \begin{bmatrix}
        0 & 1 & i  \\
        1 & 0 & -i  \\
        -i & i & 0  \\
\end{bmatrix}~\text{and}~
H(D') = \begin{bmatrix}
        0 & 0 & 0 & 0 & 0 & 0 & 0 & 0 \\
        0 & 0 & 0 & 0 & 0 & 0 & 0 & 0 \\
        0 & 0 &0 & 0 & 0 & 1 & 1 & i \\
        0 & 0 &0 & 0 & 0 & 1 & 1 & i \\
        0 & 0 &0 & 0 & 0 & 1 & 1 & i \\
        0 & 0 &1 & 1 & 1 & 0 & 0 & -i \\
        0 & 0 &1 & 1 & 1 & 0 & 0 & -i \\
        0 & 0 &-i& -i& -i& i & i & 0 
\end{bmatrix}.
\]}
We conclude this section with two observations that will be quite obvious to the experienced graph theorist, though the ideas are in a sense key to the presented discussion. As such, their proofs are omitted. 
\begin{lemma}
Let $D$ be a digraph of order $n$ and let $t$ be an expansion vector for $D$. Then
\begin{equation}
 \text{Rank }H(TR(D)) =  \text{Rank }H(D) = \text{Rank }H(TE(D,t))
\end{equation}
\end{lemma}

\begin{lemma}\label{lemma: p positive q negative ev TE}
Let $D$ be a digraph of order $n$ and let $t$ be an expansion vector for $D$. Suppose that $D$ has $p$ positive and $q$ negative eigenvalues. Then $TE(D,t)$ and $TR(D)$ also have $p$ positive and $q$ negative eigenvalues.
\end{lemma}

\section{The negative tetrahedron}
In the present paper, we are interested in families of digraphs that contain many copies of a given substructure, which is in a sense counted by the spectrum. In this section, we will provide a thorough introduction of these families and the elementary observations upon which many of the later results are built. 

\subsection{Digraphs related to the negative triangle}
\begin{figure}[b]
\begin{center}
\begin{subfigure}[b]{0.20\textwidth}\centering
\begin{tikzpicture} 
\tikzset{vertex/.style = {shape=circle,draw,minimum size=1em}}
\tikzset{edge/.style = {->,> = latex'}}
\tikzset{symedge/.style = {-,> = latex'}}
\tikzset{edgeBendLeft/.style = {->,> = latex',bend left}}
\tikzset{edgeBendRight/.style = {->,> = latex',bend right}}


\def\n{3}
\def\rad{1}
\foreach \a in {1,2,...,\n}{
\node[vertex] (\a) at (\a*360/\n: \rad cm) {};
}
\draw[symedge] (1) to node{} (2);
\draw[edge] (1) to node{} (3);
\draw[edge] (2) to node{} (3);


\end{tikzpicture}
\caption{$D_1$}
\label{fig: thm H3 1}

\end{subfigure}
\quad
\begin{subfigure}[b]{0.20\textwidth}\centering
\begin{tikzpicture} 
\tikzset{vertex/.style = {shape=circle,draw,minimum size=1em}}
\tikzset{edge/.style = {->,> = latex'}}
\tikzset{symedge/.style = {-,> = latex'}}
\tikzset{edgeBendLeft/.style = {->,> = latex',bend left}}
\tikzset{edgeBendRight/.style = {->,> = latex',bend right}}


\def\n{3}
\def\rad{1}
\foreach \a in {1,2,...,\n}{
\node[vertex] (\a) at (\a*360/\n: \rad cm) {};
}
\draw[symedge] (1) to node{} (2);
\draw[edge] (3) to node{} (1);
\draw[edge] (3) to node{} (2);


\end{tikzpicture}
\caption{$D_2$}
\label{fig: thm H3 2}

\end{subfigure}
\quad
\begin{subfigure}[b]{0.20\textwidth}\centering
\begin{tikzpicture} 
\tikzset{vertex/.style = {shape=circle,draw,minimum size=1em}}
\tikzset{edge/.style = {->,> = latex'}}
\tikzset{symedge/.style = {-,> = latex'}}
\tikzset{edgeBendLeft/.style = {->,> = latex',bend left}}
\tikzset{edgeBendRight/.style = {->,> = latex',bend right}}


\def\n{3}
\def\rad{1}
\foreach \a in {1,2,...,\n}{
\node[vertex] (\a) at (\a*360/\n: \rad cm) {};
}
\draw[symedge] (1) to node{} (2);
\draw[symedge] (1) to node{} (3);
\draw[symedge] (3) to node{} (2);


\end{tikzpicture}
\caption{$D_3$}
\label{fig: thm H3 3}

\end{subfigure}
\quad
\begin{subfigure}[b]{0.20\textwidth}\centering
\begin{tikzpicture} 
\tikzset{vertex/.style = {shape=circle,draw,minimum size=1em}}
\tikzset{edge/.style = {->,> = latex'}}
\tikzset{symedge/.style = {-,> = latex'}}
\tikzset{edgeBendLeft/.style = {->,> = latex',bend left}}
\tikzset{edgeBendRight/.style = {->,> = latex',bend right}}


\def\n{3}
\def\rad{1}
\foreach \a in {1,2,...,\n}{
\node[vertex] (\a) at (\a*360/\n: \rad cm) {};
}
\draw[symedge] (1) to node{} (2);
\draw[edge] (1) to node{} (3);
\draw[edge] (3) to node{} (2);


\end{tikzpicture}
\caption{$D_4$}
\label{fig: thm H3 4}

\end{subfigure}

\end{center}
\caption{\label{fig: thm H3} The four triangles that conribute to $\trace{H^3}$.}
\end{figure}
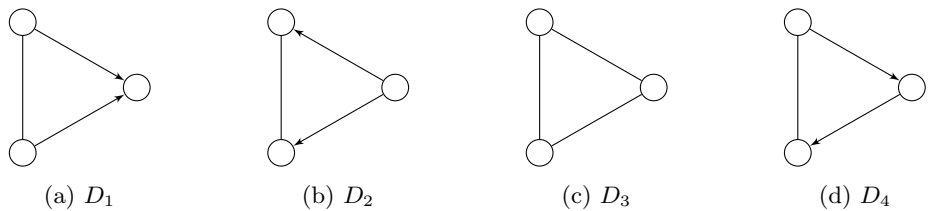

Upon investigation of properties that may be inferred from the Hermitian spectrum of a digraph, we are inspired by the following lemma by \citet{guo2015simple}, that ties in closely to a similar, well-known result for undirected graphs. 
\begin{lemma}\label{lemma: contained triangle H3}\cite{guo2015simple} Let $D$ be a digraph with Hermitian $H$. Then $\trace{H^3} = 6(x_1 + x_2 + x_3 - x_4)$, where $x_j$ denotes the number of copies of $D_j$ that occur as induced subdigraphs of $D$. The structures $D_1,\ldots,D_4$ are shown in Figure \ref{fig: thm H3}.
\end{lemma}

The main observation to take away from Lemma \ref{lemma: contained triangle H3} is that apparently, the structure $D_4$, above, is the only order 3 substructure that has a negative impact on $\trace{H^3}$, which in turn may be computed directly from the spectrum of a digraph. Thus, we may be able to identify (or even determine) digraphs that have many such substructures. In the interest of clearness, we name the following two structures, which occurred before in Figure \ref{fig: example SHDS digraph} and that are in a sense fundamental to the discussion in this paper.
\begin{definition}\label{def: neg triangle}
Figure \ref{fig: thm H3 4} is called the \textit{negative triangle} and is denoted $T^-$. 
\end{definition}
\begin{definition}
Figure \ref{fig: psi } is called the \textit{negative tetrahedron} and is denoted $K^-$.
\end{definition}
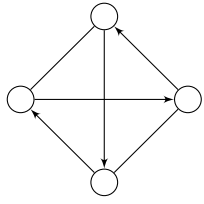
\begin{figure}[t]
     \begin{center}
     \begin{tikzpicture} [scale = 1.1]
     \tikzset{vertex/.style = {shape=circle,draw,minimum size=1em}}
     \tikzset{edge/.style = {->,> = latex'}}
     \tikzset{symedge/.style = {-,> = latex'}}
     \tikzset{edgeBendLeft/.style = {->,> = latex',bend left}}
     \tikzset{edgeBendRight/.style = {->,> = latex',bend right}}
     
     
\def\n{4}
\def\rad{1}
\foreach \a in {1,2,...,\n}{
\node[vertex] (\a) at (\a*360/\n: \rad cm) {};
}
\draw[symedge] (1) to node{} (2);
\draw[symedge] (3) to node{} (4);
\draw[edge] (1) to node{} (3);
\draw[edge] (3) to node{} (2);
\draw[edge] (2) to node{} (4);
\draw[edge] (4) to node{} (1);

     
     \end{tikzpicture}
     
     \end{center}
     \caption{\label{fig: psi } The negative tetrahedron.}
     \end{figure}

The negative tetrahedron is an interesting digraph for a number of reasons, and has come up in  the before mentioned works. One might first notice its extreme degree of structural symmetry; $K^-$ is, in fact, vertex-transitive. 
A second interesting fact is that $T^-$ and $K^-$ are exactly the two digraphs with rank more than $2$ that are antispectral\footnote{A pair of digraphs $D$ and $D'$ is said to be \textit{antispectral} to one another if $\Sigma_D = -\Sigma_{D'}$.} to a complete graph. 
$T^-$ and $K^-$, whose spectra are $\{-2,1^{[2]}\}$ and $\{-3, 1^{[3]}\}$, respectively, are antispectral to respectively $K_3$ and $K_4$.
\citet{guo2016hermitian} have shown that there are no higher rank digraphs that admit to this property. 
Lastly, it is mentioned in \cite{guo2016hermitian} that $K^-$ exhibits extreme spectral behavior, in the sense that it attains the bound $\rho(D) \leq 3\mu_1$, where $\mu_1$ is the largest eigenvalue of $D$.


In addition to $T^-$ and $K^-$, there are two more digraphs that play a prominent role throughout this paper. In the interest of structure, we include their definitions here.
\begin{definition}\label{def: Ta Tb}
The digraphs $T^-_a$ and $T^-_b$ are shown in Figures \ref{rank3reduced 1} and \ref{rank3reduced 2}, respectively. 
\end{definition}
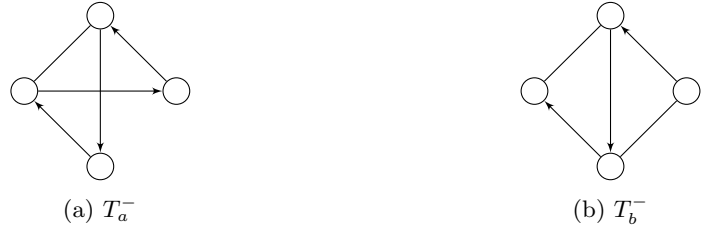
\begin{figure}[t]
     \begin{center}
     \begin{subfigure}[b]{0.45\textwidth}\centering
     \begin{tikzpicture} 
     \tikzset{vertex/.style = {shape=circle,draw,minimum size=1em}}
     \tikzset{edge/.style = {->,> = latex'}}
     \tikzset{symedge/.style = {-,> = latex'}}
     \tikzset{edgeBendLeft/.style = {->,> = latex',bend left}}
     \tikzset{edgeBendRight/.style = {->,> = latex',bend right}}
     
     
     \def\n{4}
     \def\rad{1}
     \foreach \a in {1,2,...,\n}{
     \node[vertex] (\a) at (\a*360/\n: \rad cm) {};
     }
     \draw[symedge] (1) to node{} (2);
     \draw[edge] (1) to node{} (3);
     \draw[edge] (3) to node{} (2);
     \draw[edge] (2) to node{} (4);
     \draw[edge] (4) to node{} (1);
     
     
     \end{tikzpicture}
     \caption{$T^-_a$}
     \label{rank3reduced 1}
     
     \end{subfigure}
     \begin{subfigure}[b]{0.45\textwidth}\centering
     \begin{tikzpicture} 
     \tikzset{vertex/.style = {shape=circle,draw,minimum size=1em}}
     \tikzset{edge/.style = {->,> = latex'}}
     \tikzset{symedge/.style = {-,> = latex'}}
     \tikzset{edgeBendLeft/.style = {->,> = latex',bend left}}
     \tikzset{edgeBendRight/.style = {->,> = latex',bend right}}
     
     
     \def\n{4}
     \def\rad{1}
     \foreach \a in {1,2,...,\n}{
     \node[vertex] (\a) at (\a*360/\n: \rad cm) {};
     }
     \draw[symedge] (1) to node{} (2);
     \draw[symedge] (3) to node{} (4);
     \draw[edge] (1) to node{} (3);
     \draw[edge] (3) to node{} (2);
     \draw[edge] (4) to node{} (1);
     
     
     \end{tikzpicture}
     \caption{$T^-_b$}
     \label{rank3reduced 2}
     
     \end{subfigure}
     \end{center}
     \caption{\label{rank3reduced} Illustrations for Def. \ref{def: Ta Tb}. Recall that the vertices are ordered such that '$1$' is assigned to the northmost node, then proceeding in counterclockwise order.}
     \end{figure}
We note that both digraphs are reduced and have rank $3$. 
Furthermore, we observe that $T^-_a$ and $T^-_b$ are closely related to $T^-$, from a spectral point of view. 
In fact, if one expands a single vertex of $T^-$ once (to obtain, say, $T^-_{1,1,2}$), then $T^-_{1,1,2},~ T^-_a,~ T^-_b$ are all cospectral and switching equivalent. 
More detail concerning this relation is provided at the end of Section \ref{subsec: explicit spectra}.

\subsection{The first families of SHDS digraphs}\label{subsec: first SHDS}
Using only the tools we have thus far, we are already able to construct some infinite SHDS families. 
Probably the first, most trivial SHDS digraph that comes to mind is the empty graph of order $p$, denoted $O_p$. Indeed, using that $2|E(\Gamma(O_p))| = \trace{H(O_p)^2} = 0$, the all-zero spectrum certainly determines the empty graph.

It is easy to check that $T^-$ is the smallest non-empty digraph that is strongly determined by its Hermitian spectrum. In fact, it is a simple exercise to show the following result,\footnote{A similar result occurs as Prop. 5.2 in \cite{guo2016hermitian}, which concerns digraphs antispectral to $K_n$. 
While the obtained collection of digraphs is almost identical, the authors chose to include a proof as the claim requirements and the argument are signifficantly different.} that signifies the essential role $T^-$ and $K^-$ play in the proposed discussion. 
As a first step towards a less trivial infinite family of SHDS digraphs, we classify all digraphs with largest eigenvalue $1$.  
\begin{lemma}\label{lemma: disjoint union K_n}
Let $D$ be a connected digraph with largest eigenvalue $1$. Then $D$ is either $K_2, K_2', T^-,$ or $K^-$, where $K_2'$ is the oriented $K_2$. 
\end{lemma}
\begin{proof}
By interlacing, it follows that there is no $U\in V(D)$ such that $\Gamma(D[U]) = P_3,$ since every digraph whose underlying graph is $P_3$ has $\mu_1=\sqrt{2},$ where $\mu_1$ is the largest eigenvalue. Hence, $\Gamma(D)=K_n$, for $n\in \mathbb{N}$. If $n=2$, both digraphs $D$ with $\Gamma(D) = K_2$ have $\mu_1=1$, and are therefore valid options. Note that $n=3$, $\mu_1=1$ only if $D=T^-$; other potential digraphs of order $3$ have an eigenvalue $\sqrt{3}$ or $2$.
Moreover, one should also observe that if $n>3$, each order-$3$ induced subdigraph of $D$ should again be $T^-$, or else the claim is false by interlacing. Thus, if $n=4$, $D=K^-$. 
Moreover, as before, if $n>4$, any order $4$ induced subdigraph of $D$ must be $K^-$. 

Finally, suppose that $n\geq 5$, and suppose that $D[\{1,2,3,4\}]= D[\{1,2,3,5\}] =K^-$. (Note that we are not imposing any extra assumptions; if these induced subdigraphs are not $K^-$, $D$ certainly has $\mu_1>1$.)
Then we have $D[\{1,4,5\}]\not= T^-$, as $H(D)_{1,4} = H(D)_{1,5}$. It follows that $\mu_1>1$, by which we have a contradiction and thus $n\leq 4$.
\end{proof}
By the result above, we are able to draw an interesting conclusion with regard to the spectral determination of the class of digraphs with largest eigenvalue $1$.
\begin{corollary} 
Let $D$ be a digraph with largest eigenvalue $1$. Then $D$ is WHDS. 
\end{corollary}
\begin{proof}
By Lemma \ref{lemma: disjoint union K_n}, every connected component must be $K_2, K_2', T^-,$ or $K^-$. It should be clear that every eigenvalue $-3$ belongs to a copy of $K^-$, every eigenvalue $-2$ belongs to a $T^-$, and every eigenvalue $0$ belongs to an isolated node.
Likewise, every eigenvalue $-1$ belongs to $K_2$ or $K_2'$. Since $K_2$ is clearly switching equivalent to $K_2',$ the claim follows.
\end{proof}
Moreover, if one excludes $K_2$ and $K_2'$, such a digraph is SHDS, by the same argument.
\begin{corollary}
Let $D$ be a digraph with largest eigenvalue $1$ and no eigenvalues $-1$. Then $D$ is SHDS. 
\end{corollary}
Thus, we have obtained an infinite family of SHDS digraphs, that may consist of arbitrarily many disjoint copies of $T^-$ and $K^-$, as well as isolated vertices. 

\subsection{The spectra of expansions of $K^-$ and related digraphs} \label{subsec: explicit spectra}

\begin{figure}[h!]
     \begin{center}
     \begin{tikzpicture} [scale = 0.6]
     \tikzset{vertex/.style = {shape=circle,draw,minimum size=1em,fill=black}}
     
     \tikzset{cloud/.style = {shape=circle,draw,dashed,minimum size=4.8em}}
     \tikzset{edge/.style = {->,> = latex'}}
     \tikzset{symedge/.style = {-,> = latex'}}
     \tikzset{edgeBendLeft/.style = {->,> = latex',bend left}}
     \tikzset{edgeBendRight/.style = {->,> = latex',bend right}}
     
     
     \node[cloud] (0) at (-4, 2) {};
     \node[cloud] (2) at (0, 0) {};
     \node[cloud] (3) at (4, 0) {};
     \node[cloud] (4) at (4, 4) {};
     \node[cloud] (1) at (0, 4) {};
     \draw[edge] (1) to node{} (3);
     \draw[edge] (4) to node{} (1);
     \draw[symedge] (1) to node{} (2);
     
     \draw[symedge] (3) to node{} (4);
     \draw[edge] (3) to node{} (2);
     \draw[edge] (2) to node{} (4);
     
     \node[vertex] (01) at (-4, 1.5) {};
     \node[vertex] (02) at (-4, 2.5) {};
     
     \node[vertex] (11) at (0, 0) {};
     \node[vertex] (12) at (0.7, 0) {};
     \node[vertex] (13) at (0, 0.7) {};
     \node[vertex] (14) at (-0.5, -0.5) {};
     
     \node[vertex] (21) at (3.5, 0) {};
     \node[vertex] (22) at (4.5, 0) {};
     
     \node[vertex] (31) at (0, 4) {};
     \node[vertex] (32) at (0.7, 4) {};
     \node[vertex] (31) at (-0.7, 4) {};
     \node[vertex] (31) at (0, 3.3) {};
     \node[vertex] (32) at (0, 4.7) {};
     
     \node[vertex] (41) at (4, 4.7) {};
     \node[vertex] (41) at (3.5, 3.5) {};
     \node[vertex] (42) at (4.7, 4) {};
     
     \draw (-2.5,-0.5) node {$t_2 = 4$};
     \draw (-2.5,4.5) node {$t_1 = 5$};
     \draw (6.5,4) node {$t_4 = 3$};
     \draw (6.5,0) node {$t_3 = 2$};
     \draw (-6.5,2) node {$t_0 = 2$};
     
     
     \end{tikzpicture}
     
     \end{center}
     \caption{
     \label{fig: TE(Psi)}
     A digraph obtained as 
     $TE\left(K^-,[2 ~ 5 ~ 4 ~ 2 ~ 3]
    \right)$. Here, the dashed circles indicate clusters of twins and an edge between two clusters is used to draw all edges of that type between the members of said clusters.}
     \end{figure}
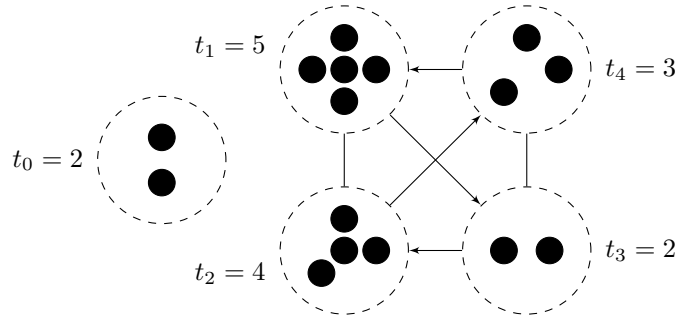

In the discussion leading up to our main theorem, we will be interested in twin expansions (recall Def. \ref{def: twin expansion}) of $K^-$, in particular. 
The following Lemma is added for completeness; its correctness should be evident from a brief look at Figure \ref{fig: TE(Psi)}. 
\begin{lemma} \label{lemma: count n m k}Let $t = \begin{bmatrix}t_0 & t_1 & t_2 & t_3 & t_4\end{bmatrix}$ and let $D_t = TE(K^-,t)$. Then $D_t$ contains 
$n = \sum_{i=0}^4 t_i$ vertices,
$m = \sum_{1\leq i < j \leq 4} t_i t_j$ edges, and 
$k = \sum_{1\leq i < j < l\leq 4} t_i t_j t_l$ copies of $T^-$.
\end{lemma}
Note that a permutation of the coefficients $t_1,\ldots,t_4$ would not necessarily yield an isomorphic digraph, as is illustrated in Example \ref{ex2}, (Sec. \ref{sec: large}) while the spectrum is invariant under such a permutation, as we will see shortly. 

As we set out to show that particular twin expansions of $K^-$ are SHDS, it seems fitting to include the explicit spectra of this interesting family of diraphs.  
In the below, we write their characteristic polynomials explicitly by employing Lemma \ref{lemma: equitable partition quotient}. While one could have used a before-mentioned result from \cite{liu2015hermitian} that counts elementary subdigraphs to obtain the coefficients in \eqref{eq: charpoly TE Psi}, the authors found the approach below to be significantly more comprehensible. 

\begin{lemma}\label{lemma: charpoly TE Psi}
Let $t=\begin{bmatrix}t_0 & t_1 & t_2 & t_3 & t_4\end{bmatrix}$ be an expansion vector and let $n = \sum_{i=0}^4 t_i$. Then $D = TE(K^-,t)$ has characteristic polynomial 
\begin{equation}\label{eq: charpoly TE Psi}
    \chi_D(\mu) = \mu^{n-4}\left(\mu^4 - \left(\sum_{1\leq i < j \leq 4} t_it_j\right)\mu^2 +2\left(\sum_{1\leq i < j < k \leq 4} t_it_jt_k\right)\mu -3\prod_{i=1}^4 t_i\right).
\end{equation}
\end{lemma}
\begin{proof}
By construction, we may write the Hermitian of $D$ as the block matrix
\begin{equation}\label{eq: partition H}
H(TE(K^-,t)) = \begin{bmatrix}0 & 0 \\ 0 & M \end{bmatrix},~\text{where }~
M = \begin{bmatrix}    0 & J & iJ & -iJ \\
                       J & 0 & -iJ & iJ \\
                      -iJ & iJ & 0 & J \\
                       iJ & -iJ & J & 0 \end{bmatrix}, \end{equation}
where the diagonal blocks have sizes $t_0\times t_0, \ldots , t_4\times t_4$, respectively.
Note that all of the blocks in \eqref{eq: partition H} are constant, and thus \eqref{eq: partition H} is a so-called equitable partition.
We may then write the $4\times 4$ quotient matrix \citep{haemers1995, guo2016hermitian} $B$ as 
\[B = \begin{bmatrix}   0 & t_2 & it_3 & -it_4\\
                       t_1 & 0 & -it_3 & it_4\\
                       -it_1 & it_2 & 0 & t_4\\
                       it_1 & -it_2 & t_3 & 0
                      \end{bmatrix}.\]
One may compute $\det(\mu I - B)$ to find
\[\chi_B(\mu) = \mu^4 - \left(\sum_{1\leq i < j \leq 4} t_it_j\right)\mu^2 +2\left(\sum_{1\leq i < j < k \leq 4} t_it_jt_k\right)\mu -3\prod_{i=1}^4 t_i.\]
Now, observe that $\rank{B}=4$. Hence, $\chi_B(\mu)$ has four nonzero roots, which are the (not necessarily distinct) eigenvalues $\lambda_1,\dots,\lambda_4$ of $B$. Since \eqref{eq: partition H} is an equitable partition, each of the $\lambda_j$ also occur as an eigenvalue of $H(TE(K^-,t))$. Moreover, since by construction $\rank{H(TE(K^-,t))} = \rank{B} = 4$, it is clear that  we have
\begin{align*}
    \chi_{D}(\mu) &= \mu^{n - 4} (\mu-\lambda_1)(\mu-\lambda_2)(\mu-\lambda_3)(\mu-\lambda_4)\\
    &= \mu^{n - 4} \chi_B(\mu).
\end{align*}
\end{proof}
Moreover, by plugging in $t$, one may readily show the following results.

\begin{corollary}\label{cor: explicit spectra}
Let $t_1,t_2,t_3\in\mathbb{N}$ and $t_0 \in \mathbb{N}_0$. For the following special cases of expansion vector $t$, we may write the spectrum of $D = TE(K^-,t)$ explicitly. In the below, $n$ is the sum of the entries of $t$. 
\begin{enumerate}[label=(\roman*)]
    \item If $t = \begin{bmatrix} t_0 & t_1 & t_1 & t_1 & t_1\end{bmatrix}$ then
\[\Sigma_{D} = \left\{-3t_1,t_1^{[3]},0^{[n - 4]} \right\}.\]
    \item If $t = \begin{bmatrix}t_0 &  t_1 & t_1 & t_1 & t_2\end{bmatrix}$ then
\[\Sigma_{D} = \left\{-t_1-\sqrt{3t_1t_2 + t_1^2},~~-t_1 + \sqrt{3t_1t_2+t_1^2},~~t_1^{[2]},~~0^{[n - 4]} \right\}.\]
    \item If $t = \begin{bmatrix} t_0 &t_1 & t_1 & t_2 & t_2\end{bmatrix}$ then
    \[\Sigma_D = \left\{ t_1, t_2, \frac{1}{2}\left(- t_1-t_2 + \sqrt{t_1^2 + 14t_1t_2 + t_2^2} \right),  \frac{1}{2}\left(-t_1-t_2 -\sqrt{t_1^2 + 14t_1t_2 + t_2^2} \right), 0^{[n - 4]} \right\}\]
\end{enumerate}
\end{corollary}
\begin{proof} 
Follows directly by plugging $t$ into Lemma \ref{lemma: charpoly TE Psi}.
\end{proof}
\begin{corollary}
Let $t_1,t_2,t_3\in\mathbb{N}$ and $t_0 \in \mathbb{N}_0$. If $t = \begin{bmatrix}t_0 &  t_1 & t_1 & t_2 & t_3\end{bmatrix}$ then $TE(K^-, t)$ has an eigenvalue $t_1$.
\end{corollary}
\begin{proof}
Plug in \eqref{eq: charpoly TE Psi} to find 
\begin{align*}
    \chi_D(\mu) &= \mu^{n - 4}(\mu-t_1)(\mu^3 + t_1 \mu^2 - ( 2t_1 t_2  + 2 t_1 t_3 + t_2t_3) \mu  + 3 t_1 t_2 t_3),
\end{align*}
which clearly has a root at $\mu=t_1.$
\end{proof}
Now, one would like to conclude that the reverse is also true; that from the occurrence of an integer eigenvalue $\mu_j$ it follows that an expansion vector contains $\mu_j$ twice. This is in general not true, as shown by the following example.
\example{ex: integer ev}{
Suppose that $t = \begin{bmatrix}0& 1& 2& 6& 9\end{bmatrix}$. Then, by Lemma \ref{lemma: charpoly TE Psi}, $D=TE(K^-,t)$ has characteristic polynomial
\[\chi_D(\mu) = \mu^{14}\left(\mu^4 -101\mu^2 + 384\mu - 324\right) = \mu^{14}(\mu-3)(\mu^3 + 3\mu^2 -92\mu + 108),\]
and thus $D$ has an eigenvalue $3$, while none of the $t_i$ equals $3.$ 
}
That said, if an integer eigenvalue $\mu_j$ occurs at least twice, then we are able to conclude the reverse, as we will discuss in the proof of Theorem \ref{thm: main theorem}.\\ 

We conclude this section with some brief notes regarding the spectral similarity of $T^-$, $T^-_a$, and $T^-_b$, and their respective expanded versions. 
As before, we are able to compute their spectra explicitly by employing the quotient matrix.
 \begin{lemma}
Let $t = \begin{bmatrix} t_0 & t_1 & t_2 & t_3 \end{bmatrix}$ be an expansion vector and let $n = \sum_i t_i$. Then $D = TE(T^-, t)$ has characteristic polynomial
\[\chi_D(\mu) = \mu^{n-3}\left(\mu^3 - (t_1t_2 + t_1t_3+ t_2t_3)\mu + 2t_1t_2t_3\right).\]
\end{lemma}
\begin{lemma}
Let $t = \begin{bmatrix} t_0 & t_1 & t_2 & t_3 & t_4 \end{bmatrix}$ and be an expansion vector. Let $D'_a = TE(T^-_a,t)$ and $D'_b = TE(T^-_b,t)$. 
Then 
\[\chi_{D'_a}(\mu) = \mu^{n-3}\left(\mu^3 - (t_1t_2 + t_1(t_3+t_4) + t_2(t_3+t_4))\mu + 2t_1t_2(t_3+t_4)\right)\]
and
\[\chi_{D'_b}(\mu) = \mu^{n-3}\left(\mu^3 - (t_1 t_3 + t_1(t_2+t_4) + t_3(t_2+t_4))\mu + 2t_1t_3(t_2+t_4)\right).\]
\end{lemma}
Thus, it follows that $TE(T^-,\begin{bmatrix}t_0 & t_1 & t_2 & (t_3 + t_4) \end{bmatrix})$,
$TE(T^-_a,\begin{bmatrix}t_0 & t_1 & t_2 & t_3 & t_4\end{bmatrix})$ and 
$TE(T^-_b,\begin{bmatrix}t_0 & t_1 & t_3 & t_2 & t_4\end{bmatrix})$ are all cospectral.
Lastly, note that these digraphs are also all pairwise switching equivalent.



\section{Classification of digraphs with one negative eigenvalue}\label{sec: small}
In order to construct the desired infinite family of SHDS digraphs, we first set out to find its smallest members. 
It turns out that the members of the family in which we are interested share the interesting property of only having a single negative eigenvalue; a property that is satisfied by very few reduced digraphs. 
We note the following useful obseravation with regard to such digraphs.
\begin{lemma}\label{lemma: 1 neg ev connected}
Let $D$ be reduced with exactly one negative eigenvalue. Then $D$ is connected.
\end{lemma}
\begin{proof}
Note that the spectrum of any connected component of order at least $2$ contains at least one negative eigenvalue, since the sum of the eigenvalues of a Hermitian matrix must sum up to its trace, which is zero for the Hermitian adjacency matrix of a digraph without loops. 
Moreover, recall that if $D$ consists of two disjoint, connected components $D_1$ and $D_2$, then $\Sigma_D = \Sigma_{D_1}\cup \Sigma_{D_2},$ and thus $\Sigma_D$ contains at least two negative elements.
Lastly, note that no isolated vertices are allowed by definition of reducedness. 
\end{proof}
We also impose a minor assumption on the rank of the considered digraphs, in order to exclude cases that are in a sense trivial. 
Specifically, we require digraphs to have rank larger than $2$. 
Recall that there are no  digraphs with rank less than $2$ besides the empty graph, and that any nonempty rank 2 digraph trivially has precisely one positive and one negative eigenvalue, by the observation above.
However, no such digraph is interesting for the present paper, as any rank $2$ digraph is cospectral to its underlying graph \cite{mohar2016hermitian} and such digraphs are in general not WDHS.\footnote{By employing twin expansion on e.g., $K_2$ and $\overrightarrow{C}_3$, one easily finds cospectral classes whose members have at least two distinct underlying graphs.}
As such, we exclude this class of digraphs; the interested reader is referred to \cite{mohar2016hermitian}, in which this class is researched in considerable detail. 

If one requires the considered digraphs to have rank larger than $2$ in addition to being reduced, one finds just four digraphs.  
The main result of this section, which is the following theorem, shows exactly that.

\begin{theorem}\label{thm: all reduced digraphs with 1 neg ev}
Let $D$ be a reduced digraph of rank larger than $2$ and with exactly one negative eigenvalue. 
Then $D$ is one of $T^-$, $T^-_a$, $T^-_b$, or $K^-$. 
\end{theorem}

In order to prove Theorem \ref{thm: all reduced digraphs with 1 neg ev}, we first show several intermediate results.
First, we will provide a few crucial observations that are somewhat obvious, but that are added for the sake of completeness.
In Proposition \ref{prop: reduced order 4 contains neg triangle class}, we will see that there are exactly three reduced digraphs on four vertices that have the required single negative eigenvalue. 
The remainder of the section is concerned with bounding the size of a reduced digraph with exactly one negative eigenvalue. In particular, we will find that such a digraph may contain at most four vertices; the correctness of Theorem \ref{thm: all reduced digraphs with 1 neg ev} then follows.

As was mentioned before, the negative triangle $T^-$ plays an essential role throughout. The simple, but useful fact that any digraph of sufficient rank must contain such triangles if it has a single negative eigenvalue, is shown below. 

\begin{lemma}\label{lemma: 1 neg ev implies D contains T-}
Let $D$ be a digraph with rank larger than $2$. If $D$ has a single negative eigenvalue, then $D[U] = T^-$ for some $U\subseteq V(D).$ 
\end{lemma}
\begin{proof}
Let $\lambda_1,\lambda_2,\ldots,\lambda_p$ be the $p$ positive eigenvalues of $D$ and let $\lambda_n$ be the negative eigenvalue. We have $\lambda_n = -\sum_j^p \lambda_j$, so
$\trace{H(D)^3}=\sum_{j}^p\lambda_j^3 -\left(\sum_{j}^p\lambda_j\right)^3<0$ and thus $D$ contains at least one negative triangle by Lemma \ref{lemma: contained triangle H3}. 
\end{proof}

\begin{corollary}\label{cor: n=3, rank=3 => D = T-}
If $D$ is a digraph with order $3$, rank larger than $2$, and exactly one negative eigenvalue, then $D = T^-$.
\end{corollary}

We are now ready to show the first major result, in which we obtain the collection of order four digraphs that meets our requirements. We would remark here that, while the collection of all order four digraphs is small enough to simply apply full enumeration by computer, the authors opt for a constructive argument that may be used in similar fashion when the order is increased. 

\begin{proposition}\label{prop: reduced order 4 contains neg triangle class}
Let $D$ be a reduced digraph with order $4$, rank larger than $2$, and exactly one negative eigenvalue
Then $D$ is one of $T^-_a$, $T^-_b$, or $K^-$. 
\end{proposition}
     
\begin{proof}
     Let $D$ be a digraph of order $4$, with exactly one negative eigenvalue, and rank larger than $2$. By Lemma \ref{lemma: 1 neg ev implies D contains T-}, $D$ contains $T^-$. Hence, we may write $H(D)$ as
     \[H(D) = \begin{bmatrix} 0 & 1 & i & -i\cdot \bar{z}_1\\
                              1 & 0 &-i & i\cdot \bar{z}_2\\
                             -i & i & 0 & \bar{z}_3\\
                              i\cdot z_1 & -i \cdot z_2 & z_3 & 0\end{bmatrix},\]
    where $z = [z_1~~z_2~~z_3]\not= 0 $, $z_j \in\{0,\pm 1,\pm i\}$ and $z_1 \not=i,$ $z_2 \not=-i,$ $z_3 \not=-1.$ 
    Note that the variable entries of $H(D)$ are put in this form to make \eqref{eq: det HS} symmetric. One then readily obtains that 
    \begin{equation}\label{eq: det HS}\det H(D) = \sum_{j=1}^3 ~|z_j|~ -2 \sum_{1\leq i<j\leq 3} \re{z_i\bar{z}_j}.\end{equation}
    Now, we make the following observation. Since $D$ contains $T^-$, it has at least two positive eigenvalues, by interlacing.
    Moreover, from $\det H(D)  = \prod_{j}\lambda_j$ it follows that $\det H(D) >0$ if and only if $D$ has an even number of negative eigenvalues. 
    Hence, $D$ satisfies the requirements of the claim if and only if the corresponding $z$ is such that $\det H(D) \leq 0$. 
    Note that if exactly one $z_j$ is nonzero, one may plug in \eqref{eq: det HS} to obtain $\det H(D) =1,$ and thus $H(D)$ has more than one negative eigenvalue.
    Therefore, at least two elements of $z$ must be nonzero. 
    
    Suppose that two elements of $z$ are nonzero. Then $\det H(D)  \geq 0$ (by \eqref{eq: det HS}) and thus we are only interested in the case that $\det H(D)  = 0.$ 
    Suppose that $z_3 = 0.$ Then $\det H(D)  = 0 \iff |z_1| + |z_2| = 2\re{z_1\bar{z}_2} \iff z_1\bar{z}_2=1 \iff z_1 = z_2.$ Hence, either $z = \begin{bmatrix} 1 & 1 & 0\end{bmatrix}$ or $z = \begin{bmatrix} -1 & -1 & 0\end{bmatrix}$. 
    Similarly, if $z_2 = 0$ then $z = \begin{bmatrix} 1 & 0 & 1\end{bmatrix}$ or $z = \begin{bmatrix} -i & 0 & -i\end{bmatrix}$ and 
    if $z_1 = 0$ then $z = \begin{bmatrix} 0 & 1 & 1\end{bmatrix}$ or $z = \begin{bmatrix} 0 & i & i\end{bmatrix}$. It is readily verified that out of these six possible $z$, three correspond to a digraph that contains a twin and therefore do not meet the requirements of the claim; the remaining three $z$ correspond to either $T^-_a$ or $T^-_b$. 
    
    Lastly, suppose that no element of $z$ is zero. It is easily observed from \eqref{eq: det HS} that $\det H(D) \not=0,$ since $\re{z_i\bar{z}_j}\in \mathbb{Z}$ $\forall i, j$ and $3=2(m_1 + m_2 + m_3)$ has no solution for $m_1, m_2, m_3\in\mathbb{Z}$. Thus, any $z$ that meets the requirements of the claim has $\det H(D) <0.$ Hence, $\re{z_i\bar{z}_j} \geq 0$ $\forall i,j$, with at most one pair $(i,j)$ such that $\re{z_i\bar{z}_j} = 0$. W.l.o.g., assume that $z_1\bar{z}_2 = z_1\bar{z}_3 = 1$. Then $z_1 = z_2$ and $z_1 = z_3$, and thus $z = \begin{bmatrix}1 & 1 & 1 \end{bmatrix},$ which corresponds to exactly $K^-$. 
\end{proof}

In order for us to prove Theorem \ref{thm: all reduced digraphs with 1 neg ev}, we should consider the reduced digraphs of larger order. 
As we will shortly show, we find that no digraph of order larger than $4$ satisfies both required properties, i.e., being reduced and having exactly one negative eigenvalue. 
In the interest of structure, the discussion to support this claim is split up into Lemmas \ref{lemma: if not contain Ta Tb Psi and n>=4, then twin}, \ref{lemma: if n=5 and 1 neg eigenvalue, then twin} and \ref{lemma: if n>5 and 1 neg eigenvalue, then twin}. 


The approach below is, for the most part, based on the idea of taking some small substructure that is certainly contained in any of the candidates, and attempting to build a digraph that meets all requirements by adding vertices and edges to it. 
In particular, all of the digraphs we encounter contain at least one copy of $T^-$. 
Moreover, from the results at the top of this section, we know that there are scarcely any ways to extend $T^-$ with vertices and edges without invalidating assumptions. 
Using these facts, we will show in Lemmas \ref{lemma: if not contain Ta Tb Psi and n>=4, then twin}, \ref{lemma: if n=5 and 1 neg eigenvalue, then twin} and \ref{lemma: if n>5 and 1 neg eigenvalue, then twin} that the order of the digraphs that have so far been shown to satisfy our requirements cannot be extended without introducing a twin vertex.

\begin{lemma}\label{lemma: if not contain Ta Tb Psi and n>=4, then twin}
Suppose that $D$ is a digraph of order $n\geq 4$, rank larger than $2$, that does not contain $T^-_a$, $T^-_b$,  or $K^-$. Moreover, suppose that $D$ has exactly one negative eigenvalue. Then $D$ is not reduced.
\end{lemma}
\begin{proof}
By contradiction. We will show, through combinatorial reasoning, that a digraph that admits to the assumptions in the claim must contain twins.
This reasoning is illustrated with an example in Figure \ref{fig: proof lemma not contain Ta Tb Psi}.
We note that while the exact nature of the edges in Figure \ref{fig: proof lemma not contain Ta Tb Psi} may differ depending on $u$, the reasoning below remains valid.\footnote{In fact, one could even disregard the arc orientations in Figure \ref{fig: proof lemma not contain Ta Tb Psi} altogether and apply the reasoning in the proof of Lemma \ref{lemma: if not contain Ta Tb Psi and n>=4, then twin} to the underlying graph of the considered digraph, without compromising the proof. The authors have opted for an example that best illustrates the situation at hand.}

Suppose that $D$ is reduced. By Lemma \ref{lemma: 1 neg ev implies D contains T-}, $D$ contains $T^-$ as an induced subdigraph and $D$ is connected by Lemma \ref{lemma: 1 neg ev connected}.
Let $U\subset V(D)$ be such that $D[U] = T^-$.
As $n\geq 4$, there must be a vertex $v\in v(D)$ that is adjacent to $U$ in $D$. 
We now observe that $D[U\cup\{v\}]$ must contain a twin.
Indeed, if we suppose that $D[U\cup \{v\}]$ is reduced, then by Proposition \ref{prop: reduced order 4 contains neg triangle class}, $D[U\cup \{v\}]$ is either $T^-_a$, $T^-_b$ or $K^-$. 
However, since none of these digraphs are allowed to be contained as induced subdigraphs, by the requirements of the claim, we have a contradiction. 
Thus, $D[U\cup\{v\}]$ is not reduced, and $v$ is the twin of a vertex $u\in U$ in $D[U\cup\{v\}]$.
(At this point, $D[U\cup\{v\}]$ may look like Figure \ref{fig: proof lemma not contain Ta Tb Psi 1}.)

But since $D$ is assumed to be reduced, there must be some vertex that distinguishes $u$ from $v$ in $D$. 
Let us call this vertex $w$ and assume without loss of generality that $w$ is adjacent to $u$. 
Then, consider $D[U \cup \{w\}]$ and use the same argument as above to obtain that $w$ must be the twin of some vertex $x \in U\setminus \{u\}$ (since $u$ and $w$ are adjacent) in $D[U\cup\{w\}].$
Label the final unlabeled vertex in $U$ with $y$; we then obtain Figure \ref{fig: proof lemma not contain Ta Tb Psi 2}.
Note that, again, by the same argument, $D[\{v,w,x,y\}]$ contains $T^-$ and should thus contain a twin. 
But as $w$ is not adjacent to $x$ (since $w$ is the twin of $x$ in $D[\{u,w,x,y\}]$), $w$ must be adjacent to $v$ and implicitly $w$ is the twin of $x$ in $D[\{v,w,x,y\}]$ as well, (Figure \ref{fig: proof lemma not contain Ta Tb Psi 3}) and hence also in $D[\{u,v,w,x,y\}]$.
Finally, if we consider the full structure (Figure \ref{fig: proof lemma not contain Ta Tb Psi 4}) it then follows that $w$ does not distinguish between $u$ and $v$, which is a contradiction. 
\end{proof}
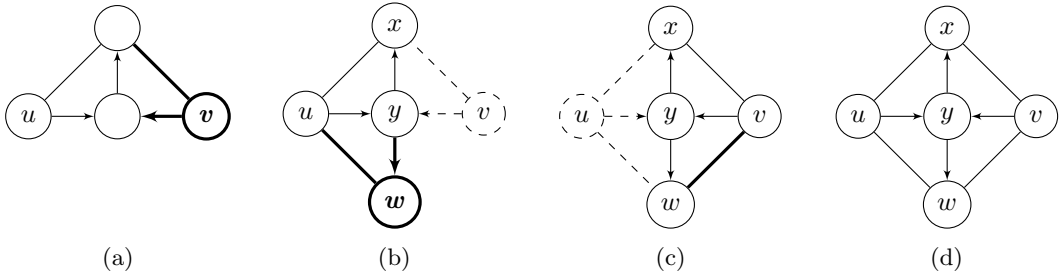
\begin{figure}[h!]
     \begin{center}
     \begin{subfigure}[b]{0.24\textwidth}\centering
     \begin{tikzpicture} [scale=0.9]
     \tikzset{vertex/.style = {shape=circle,draw,minimum size=1.7em}}
     \tikzset{edge/.style = {->,> = latex'}}
     \tikzset{symedge/.style = {-,> = latex'}}
     \tikzset{edgeBendLeft/.style = {->,> = latex',bend left}}
     \tikzset{edgeBendRight/.style = {->,> = latex',bend right}}
     \tikzset{cloud/.style = {shape=circle,draw,dashed,minimum size=6.8em}}
     \tikzset{vertex_fade/.style = {shape=circle,draw,dashed,minimum size=1em}}
     \tikzset{edge_fade/.style = {->, >=latex',dashed}}
     \tikzset{symedge_fade/.style = {-,dashed}}
     
     \def\n{4}
     \def\rad{1.3}
     \node[vertex] (x) at (1*360/\n: \rad cm) {};
     \node[vertex] (u) at (2*360/\n: \rad cm) {$u$};
     \node[vertex, color=white] (w) at (3*360/\n: \rad cm) {$w$};
     \node[vertex, very thick] (v) at (4*360/\n: \rad cm) {\textbf{\textit{v}}};
     \node[vertex] (y) at (0,0) {};
     
     
     \draw[edge] (u) to node{} (y);
     \draw[edge, very thick] (v) to node{} (y);
     \draw[edge] (y) to node{} (x);
     
     \draw[symedge] (x) to node{} (u);
     \draw[symedge, very thick] (v) to node{} (x);
     
     
     \end{tikzpicture}
     \caption{}
     \label{fig: proof lemma not contain Ta Tb Psi 1}
     
     \end{subfigure} 
     \begin{subfigure}[b]{0.24\textwidth}\centering
     \begin{tikzpicture} [scale=0.9]
     \tikzset{vertex/.style = {shape=circle,draw,minimum size=1.7em}}
     \tikzset{edge/.style = {->,> = latex'}}
     \tikzset{symedge/.style = {-,> = latex'}}
     \tikzset{edgeBendLeft/.style = {->,> = latex',bend left}}
     \tikzset{edgeBendRight/.style = {->,> = latex',bend right}}
     \tikzset{cloud/.style = {shape=circle,draw,dashed,minimum size=4.8em}}
     \tikzset{vertex_fade/.style = {shape=circle,draw,dashed,minimum size=1em}}
     \tikzset{edge_fade/.style = {->, >=latex',dashed}}
     \tikzset{symedge_fade/.style = {-,dashed}}
     
     \def\n{4}
     \def\rad{1.3}
     \node[vertex] (x) at (1*360/\n: \rad cm) {$x$};
     \node[vertex] (u) at (2*360/\n: \rad cm) {$u$};
     \node[vertex, very thick] (w) at (3*360/\n: \rad cm) {\textbf{\textit{w}}};
     \node[vertex_fade] (v) at (4*360/\n: \rad cm) {$v$};
     \node[vertex] (y) at (0,0) {$y$};
     
     \draw[edge] (u) to node{} (y);
     \draw[edge_fade] (v) to node{} (y);
     \draw[edge] (y) to node{} (x);
     \draw[edge, very thick] (y) to node{} (w);
     
     \draw[symedge] (x) to node{} (u);
     \draw[symedge, very thick] (u) to node{} (w);
     \draw[symedge_fade] (v) to node{} (x);
     
     
     \end{tikzpicture}
     \caption{}
     \label{fig: proof lemma not contain Ta Tb Psi 2}
     
     \end{subfigure}
     \begin{subfigure}[b]{0.24\textwidth}\centering
     \begin{tikzpicture} [scale=0.9]
     \tikzset{vertex/.style = {shape=circle,draw,minimum size=1em}}
     \tikzset{edge/.style = {->,> = latex'}}
     \tikzset{symedge/.style = {-,> = latex'}}
     \tikzset{edgeBendLeft/.style = {->,> = latex',bend left}}
     \tikzset{edgeBendRight/.style = {->,> = latex',bend right}}
     \tikzset{cloud/.style = {shape=circle,draw,dashed,minimum size=4.8em}}
     \tikzset{vertex_fade/.style = {shape=circle,draw,dashed,minimum size=1em}}
     \tikzset{edge_fade/.style = {->, >=latex',dashed}}
     \tikzset{symedge_fade/.style = {-,dashed}}
     
     \def\n{4}
     \def\rad{1.3}
     \node[vertex] (x) at (1*360/\n: \rad cm) {$x$};
     \node[vertex_fade] (u) at (2*360/\n: \rad cm) {$u$};
     \node[vertex] (w) at (3*360/\n: \rad cm) {$w$};
     \node[vertex] (v) at (4*360/\n: \rad cm) {$v$};
     \node[vertex] (y) at (0,0) {$y$};
     
     \draw[edge_fade] (u) to node{} (y);
     \draw[edge] (v) to node{} (y);
     \draw[edge] (y) to node{} (x);
     \draw[edge] (y) to node{} (w);
     
     \draw[symedge_fade] (x) to node{} (u);
     \draw[symedge_fade] (u) to node{} (w);
     \draw[symedge, very thick] (w) to node{} (v);
     \draw[symedge] (v) to node{} (x);
     
     
     \end{tikzpicture}
     \caption{}
     \label{fig: proof lemma not contain Ta Tb Psi 3}
     
     \end{subfigure}
     \begin{subfigure}[b]{0.24\textwidth}\centering
     \begin{tikzpicture} [scale=0.9]
     \tikzset{vertex/.style = {shape=circle,draw,minimum size=1em}}
     \tikzset{edge/.style = {->,> = latex'}}
     \tikzset{symedge/.style = {-,> = latex'}}
     \tikzset{edgeBendLeft/.style = {->,> = latex',bend left}}
     \tikzset{edgeBendRight/.style = {->,> = latex',bend right}}
     \tikzset{cloud/.style = {shape=circle,draw,dashed,minimum size=4.8em}}
     \tikzset{vertex_fade/.style = {shape=circle,draw,dashed,minimum size=1em}}
     \tikzset{edge_fade/.style = {->, >=latex',dashed}}
     \tikzset{symedge_fade/.style = {-,dashed}}
     
     \def\n{4}
     \def\rad{1.3}
     \node[vertex] (x) at (1*360/\n: \rad cm) {$x$};
     \node[vertex] (u) at (2*360/\n: \rad cm) {$u$};
     \node[vertex] (w) at (3*360/\n: \rad cm) {$w$};
     \node[vertex] (v) at (4*360/\n: \rad cm) {$v$};
     \node[vertex] (y) at (0,0) {$y$};
     
     \draw[edge] (u) to node{} (y);
     \draw[edge] (v) to node{} (y);
     \draw[edge] (y) to node{} (x);
     \draw[edge] (y) to node{} (w);
     
     \draw[symedge] (x) to node{} (u);
     \draw[symedge] (u) to node{} (w);
     \draw[symedge] (w) to node{} (v);
     \draw[symedge] (v) to node{} (x);
     
     
     \end{tikzpicture}
     \caption{}
     \label{fig: proof lemma not contain Ta Tb Psi 4}
     
     \end{subfigure}
     \end{center}
     \caption{\label{fig: proof lemma not contain Ta Tb Psi}Illustration for the proof of Lemma \ref{lemma: if not contain Ta Tb Psi and n>=4, then twin}. Bold elements indicate the most recent additions to the structure; dashed elements are temporarily ignored.}
     \end{figure}
     
Let us now consider digraphs that do contain $T^-_a$, $T^-_b$, or $K^-$.

\begin{lemma}\label{lemma: if n=5 and 1 neg eigenvalue, then twin}
Suppose that $D$ is a digraph of order $n=5$, that contains $T^-_a$, $T^-_b$, or $K^-$. Moreover, suppose that $D$ has exactly one negative eigenvalue. Then $D$ is not reduced.
\end{lemma}
\begin{proof}
First, suppose that $U\subset V(D)$ is such that $D[U] \cong K^-$.
Let $v$ be the fifth vertex in $V(D)$. 
We may assume that $v$ is not an isolated vertex, otherwise $D$ would not be reduced. 
We make the following observations from Proposition \ref{prop: reduced order 4 contains neg triangle class}: $v$ cannot have valency 1, and if $v$ has valency $3$, then the subdigraph of $D$ induced by $v$ and its neighbors is isomorphic to $K^-$. 
It then follows that $v$ is connected  to at least two out of every three vertices in $U$.
This, in turn, implies that $v$ has valency at least $3$. 

Suppose that $v$ has valency $3$, let $u\in U$ be non-adjacent to $v$ and let $U' = U\setminus \{u\} \cup \{v\}$. Then $D[U'] \cong K^-$, and thus it follows that $u$ and $v$ are twins.
If $v$ has valency $4$, one may apply the same argument twice to obtain that $v$ should be the twin of two distinct vertices in $U$, which is impossible. 
Hence, $v$ has valency $3$ and $D$ is not reduced.

Next, suppose that $D[U] \cong T^-_a.$ Then one may write 
\setbool{@fleqn}{false}
\begin{align}
 \det H(D)=&~ \det
 \begin{bmatrix}
 0 & 1 & i & -i & -i\cdot \bar{z}_1  \\
 1 & 0 & -i & i & i\cdot \bar{z}_2  \\
 -i & i & 0 & 0 & \bar{z}_3  \\
 i & -i & 0 & 0 & \bar{z}_4  \\
 i\cdot z_1 & -i\cdot z_2 & z_3 & z_4 & 0  
 \end{bmatrix} =\det
 \begin{bmatrix}
 0 & 1 & i & 0 & \ldots  \\
 1 & 0 & -i & 0 & \ldots  \\
 -i & i & 0 & 0 & \ldots  \\
 0 & 0 & 0 & 0 & \bar{z}_3 + \bar{z}_4  \\
 \ldots & \ldots & \ldots & z_3 + z_4 & 0  
 \end{bmatrix}\nonumber  \\
 =&~ -|z_3 + z_4|^2 \det
 \begin{bmatrix}
 0 & 1 & i  \\
 1 & 0 & -i  \\
 -i & i & 0 
 \end{bmatrix}  = 2|z_3 + z_4|^2,   \label{eq: neg determinant T-a}
 \end{align}
 where $z_j\in\{0,\pm 1, \pm i\}$ ($j\in[4]$), $z_1 \not=i$, $z_2 \not= -i$ and $z_3,z_4 \not= -1.$
 Since a positive determinant implies an even (and thus larger than $1$) number of negative eigenvalues, it follows that $z_3 + z_4 = 0$. Note that the proof of Proposition \ref{prop: reduced order 4 contains neg triangle class} gives us all possible solutions for $\begin{bmatrix} z_1 & z_2 & z_3\end{bmatrix}$. We remove those with $z_3 = 1$, since $z_4 = -1$ is not allowed. Besides $z = 0$, we obtain $z = \begin{bmatrix}1 & 1 & 0 & 0 \end{bmatrix}$,  
 $\begin{bmatrix}-1 & -1 & 0 & 0 \end{bmatrix}$,
 $\begin{bmatrix}-i & 0 & -i & i \end{bmatrix}$,  
 $\begin{bmatrix}0 & i & i & -i \end{bmatrix}$. The solution $z=0$ makes the fifth vertex an isolated vertex, whereas each of the other solutions makes it a twin of one of the four other vertices. 
 
  The proof for $D[U] \cong T^-_b$ is analogous to the above and therefore omitted.
 \setbool{@fleqn}{true}
\end{proof}

Using Lemma \ref{lemma: if n=5 and 1 neg eigenvalue, then twin}, we are able to extend its claim to arbitrarily large $n$. 

\begin{lemma}\label{lemma: if n>5 and 1 neg eigenvalue, then twin}
Suppose that $D$ is a digraph of order $n>5$, that contains $T^-_a$, $T^-_b$ or $K^-$. Moreover, suppose that $D$ has exactly one negative eigenvalue. Then $D$ is not reduced.
\end{lemma}
\begin{proof}
By contradiction. Suppose that $D$ is reduced and let $U\subset V(D)$ be such that $D[U]$ is either $K^-$, $T^-_a$ or $T^-_b$. Fix some $u' \in V(D)\setminus U$ such that $U^* = U\cup\{u'\}$ induces a weakly connected subdigraph $D[U^*]$. By Lemma \ref{lemma: if n=5 and 1 neg eigenvalue, then twin}, there is a vertex in $U$, say $u$, such that $u$ and $u'$ are twins in $D[U^*].$ If we let $U' = U^*\setminus \{u\},$ then $D[U] \cong D[U'].$

Because $D$ is assumed to be reduced, there is some vertex $w\in V(D)\setminus U^*$ that distinguishes $u$ from $u'.$ Let $W = U\cup \{w\}$ and $W' = U' \cup \{w\}$ and observe that, by Lemma \ref{lemma: if n=5 and 1 neg eigenvalue, then twin}, $w$ is twin to a member of $U$ and $U'$ in $D[W]$ and $D[W']$, respectively. Now, let $v$ be a vertex in $U\setminus \{u\}$ with valency $3$ in $D[U]$. Note that regardless of the choice of $u$, such a vertex exists for each of the cases for $D[U]$.  Moreover, the relation to $v$ (in-neighbor, out-neighbor, undirected neighbor or no neighbor) is different for each of the four vertices in $D[U]$. This implies that the relation of $w$ to $v$ determines which of the vertices in $U$ is the twin of $w$ in $D[W]$. 

Suppose that the twin of $w$ in $D[W]$ is $u$. Then the twin of $w$ in $D[W']$ is $u'$, because $u$ and $u'$ have the same relation to $v$. But then $w$ does not distinguish $u$ from $u'$, since it is not adjacent to either of the two. This means that the twin of $w$ in $D[W]$ must be a member of $U\setminus \{u\} = U'\setminus \{u'\}$. Specifically, it is the same vertex as the twin of $w$ in $D[W']$, as a consequence of the unique relation to $v$. But this twin does not distinguish $u$ from $u'$ (as $D[U] \cong D[U']$), so neither does $w$, and we have our final contradiction. Thus, $D$ is not reduced. 
\end{proof}

By Lemmas \ref{lemma: if not contain Ta Tb Psi and n>=4, then twin} through \ref{lemma: if n>5 and 1 neg eigenvalue, then twin}, we now have all the necessary tools to prove Theorem \ref{thm: all reduced digraphs with 1 neg ev}.

\begin{proof}\textit{(Of Theorem \ref{thm: all reduced digraphs with 1 neg ev}.)}
Suppose $D$ is a  reduced digraph of order $n$ with rank larger than $2$ and exactly one negative eigenvalue.
Since any digraph of order at most $2$ has rank at most $2$, it follows that $n\geq 3.$ 
Next, we distinguish two cases: either $D$ contains at least one of $T^-_a$, $T^-_b$, $K^-$ or it does not.
In the former case, we have by Lemmas \ref{lemma: if n=5 and 1 neg eigenvalue, then twin} and \ref{lemma: if n>5 and 1 neg eigenvalue, then twin} that $n\leq 4,$ otherwise we would lose the reducedness of $D$.
In the latter case, the same conclusion follows unless $n\leq 3,$ by Lemma \ref{lemma: if not contain Ta Tb Psi and n>=4, then twin}.

We are thus left with two possibilities: either $n=3$ or $n=4$. If $n=3$, we have  that $D = T^-$, by Corollary \ref{cor: n=3, rank=3 => D = T-}. Finally, if $n=4,$ then by Proposition \ref{prop: reduced order 4 contains neg triangle class} it holds that $D\in\{T^-_a, T^-_b, K^-\}.$
\end{proof}

We are, in fact, able to conclude much more.
Using that the numbers of positive and negative eigenvalues do not change when twin reduction is applied, the results of Theorem \ref{thm: all reduced digraphs with 1 neg ev} extend to the underlying structure of any digraph with a single negative eigenvalue. 
This key observation is formalized in Theorem \ref{thm: TR of a digraph with 1 neg ev}. 


\begin{theorem}\label{thm: TR of a digraph with 1 neg ev}
Let $D$ be a digraph of order $n\geq 5$, rank larger than $2$, with exactly one negative eigenvalue. Then one of the following cases is true. 
\begin{enumerate}[label=(\roman*)]
\item $\rank{D} = 3$ and either $TR(D) = T^-$, $TR(D) = T^-_a$ or $TR(D) = T^-_b$,
\item $\rank{D} = 4$ and $TR(D) = K^-$. 
\end{enumerate}
\end{theorem}
\begin{proof}
Let $D' = TR(D)$. Then $D'$ is reduced and has exactly one negative eigenvalue, which by Theorem \ref{thm: all reduced digraphs with 1 neg ev} implies that $D'\in\{T^-, T^-_a, T^-_b, K^-\}.$ The claim clearly follows. 
\end{proof}

In particular, we observe that if one is given a spectrum that contains three positive, one negative, and arbitrarily many zero eigenvalues, one could say with certainty that the underying structure of the corresponding digraph is exactly $K^-$. 
In other words, this digraph is a twin expansion of $K^-.$
Inspired by this fact, the authors were convinced that many SHDS digraphs were within reach. 
Consider, for example, a spectrum of the form $$\Sigma = \left\{t_1^{[3]},0^{[4t_1 + t_0 - 4]},-3t_1\right\} \text{~for~} t_1\in\mathbb{N}, t_0\in\mathbb{N}_0.$$
It is now not hard to show that  this spectrum belongs to  $D = TE\left(K^-, \begin{bmatrix}t_0 & t_1 & t_1 & t_1 & t_1\end{bmatrix}\right)$ by considering $\sum_{\mu\in\Sigma} \mu^2$, $\sum_{\mu\in\Sigma} \mu^3$ and Theorem \ref{thm: TR of a digraph with 1 neg ev}. 
In the next section, we show a more general result, based on these principles.\\


As mentioned before, any rank $2$ digraph trivially has precisely one negative eigenvalue. 
For completeness, we recall that a digraph has rank 2 if and only if $\Gamma(TR(D))$ is either $K_2$, $P_3$, or $K_3$, where in the latter case it must additionally be required that $TR(D)$ contains an odd number of arcs.
For more detail, the interested reader is referred to the full characterization in \cite{mohar2016hermitian}.

\section{An infinite family of connected SHDS digraphs}\label{sec: large}
We have so far restricted ourselves to \textit{reduced} digraphs, as digraphs that admit to this assumption are in a sense the fundamental structure to the digraphs that may be obtained by introducing twin vertices. In this section, we will be relaxing said assumption and consider twin expansions of the digraphs from Theorem \ref{thm: all reduced digraphs with 1 neg ev}, to further inquire into the class of digraphs with exactly one negative eigenvalue. In particular, we use that there is exactly one of those digraphs with rank four, to arrive at a remarkable conclusion.

The main result of this section is Theorem \ref{thm: main theorem}, in which we show that any digraph $D = TE(K^-, t)$, where $t = \begin{bmatrix}t_0 & t_1 & t_1 & t_1 & t_2 \end{bmatrix}$ is an expansion vector, is strongly determined by its Hermitian spectrum. In other words, we obtain an infinite family of digraphs that is SHDS, which includes the first such connected infinite family. 

As was briefly touched on in the introduction, such digraphs are incredibly rare, as there is an extreme degree of similarity within the collection of Hermitian adjacency matrices of a given order, informally speaking.
Even when one just considers a digraph and its converse, which are clearly cospectral but in general not isomorphic, and hence in general not SHDS.
Indeed, one is easily convinced that any SHDS digraph is necessarily \textit{self-converse}, which is by itself an extremely rare property. 
While no formal proof has been found in the literature, one is easily convinced that the claim "the fraction of self-converse digraphs of order $n$ goes to zero as $n$ goes to infinity" should be true.
Indeed, this conclusion seems certain if one numerically evaluates the counting polynomials by \citet{harary1955number} and \citet{harary1966enumeration}, which count the number of digraphs and the number of self-converse digraphs on $n$ vertices, respectively. 
This evaluation is included in Appendix \ref{sec: fraction self-converse}. 

\begin{figure}[h!]
     \begin{center}
     \begin{subfigure}[b]{0.40\textwidth}\centering
     \begin{tikzpicture} 
     \tikzset{vertex/.style = {shape=circle,draw,minimum size=1em}}
     \tikzset{edge/.style = {->,> = latex'}}
     \tikzset{symedge/.style = {-,> = latex'}}
     \tikzset{edgeBendLeft/.style = {->,> = latex',bend left}}
     \tikzset{edgeBendRight/.style = {->,> = latex',bend right}}
     
     \def\n{4}
     \def\rad{1.3}
     \node[vertex] (11) at (0.8*360/\n: \rad cm) {};
     \node[vertex] (12) at (1.2*360/\n: \rad cm) {};
     \node[vertex] (21) at (1.8*360/\n: \rad cm) {};
     \node[vertex] (22) at (2.2*360/\n: \rad cm) {};
     \node[vertex] (3) at (3*360/\n: \rad cm) {};
     \node[vertex] (4) at (4*360/\n: \rad cm) {};
     
     \draw[symedge] (11) to node{} (21);
     \draw[symedge] (12) to node{} (21);
     \draw[symedge] (11) to node{} (22);
     \draw[symedge] (12) to node{} (22);
     \draw[edge] (3) to node{} (21);
     \draw[edge] (3) to node{} (22);
     \draw[symedge] (3) to node{} (4);
     \draw[edge] (4) to node{} (11);
     \draw[edge] (4) to node{} (12);
     
     \draw[edge] (11) to node{} (3);
     \draw[edge] (12) to node{} (3);
     
     \draw[edge] (21) to node{} (4);
     \draw[edge] (22) to node{} (4);
     
     \end{tikzpicture}
     \caption{$D$}
     \label{fix ex2 1}
     \end{subfigure}
     \begin{subfigure}[b]{0.40\textwidth}\centering
     \begin{tikzpicture} 
     \tikzset{vertex/.style = {shape=circle,draw,minimum size=1em}}
     \tikzset{edge/.style = {->,> = latex'}}
     \tikzset{symedge/.style = {-,> = latex'}}
     \tikzset{edgeBendLeft/.style = {->,> = latex',bend left}}
     \tikzset{edgeBendRight/.style = {->,> = latex',bend right}}
     
     \def\n{4}
     \def\rad{1.3}
     \node[vertex] (11) at (0.8*360/\n: \rad cm) {};
     \node[vertex] (12) at (1.2*360/\n: \rad cm) {};
     \node[vertex] (2) at (2*360/\n: \rad cm) {};
     \node[vertex] (31) at (2.8*360/\n: \rad cm) {};
     \node[vertex] (32) at (3.2*360/\n: \rad cm) {};
     \node[vertex] (4) at (4*360/\n: \rad cm) {};
     
     \draw[symedge] (11) to node{} (2);
     \draw[symedge] (12) to node{} (2);
     \draw[edge] (31) to node{} (2);
     \draw[edge] (32) to node{} (2);
     \draw[symedge] (31) to node{} (4);
     \draw[symedge] (32) to node{} (4);
     \draw[edge] (4) to node{} (11);
     \draw[edge] (4) to node{} (12);
     
     \draw[edge] (11) to node{} (31);
     \draw[edge] (12) to node{} (31);
     \draw[edge] (11) to node{} (32);
     \draw[edge] (12) to node{} (32);
     
     \draw[edge] (2) to node{} (4);
     
     \end{tikzpicture}
     \caption{$D'$}
     \label{fix ex2 2}
     \end{subfigure}
     \end{center}
     \caption{\label{fig ex2} Digraphs illustrated in Example \ref{ex2}.}
     \end{figure}
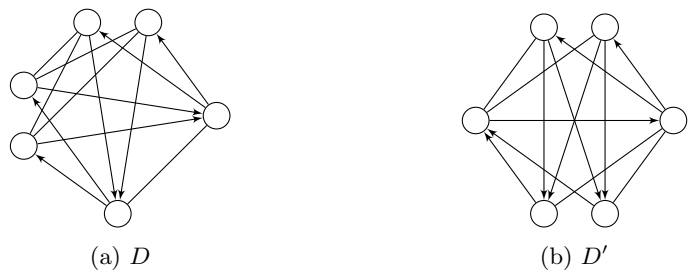

 First, we observe that, regretfully, there are still many twin expansions of $K^-$ that we may not determine uniquely from their spectra, as the following example illustrates.
\example{ex2}{
Let $t = \begin{bmatrix} 0& 2 & 2 & 1 & 1\end{bmatrix}$ and $t' = \begin{bmatrix}0& 2 & 1 & 2 & 1\end{bmatrix}.$ Then $D=TE(K^-,t)$ and $D'=TE(K^-,t')$ are $H$-cospectral by Lemma \ref{lemma: charpoly TE Psi}. However, they are clearly not isomorphic, as is visible in Figure \ref{fig ex2}: $D$ contains $5$ digons, whereas $D'$ contains only $4$.
}

Example \ref{ex2} clearly illustrates the main obstacle in this part of our quest to construct SHDS digraphs; expansion vectors that are closely related, but which do not quite yield isomorphic digraphs when used to expand $K^-$. Specifically, permutations of the same expansion vector yield cospectral, but not necessarily isomorphic, digraphs. This observation is formalized in the following lemmas. 



\begin{lemma}\label{lemma: new permute beta part 1}
Let $t_0\in\mathbb{N}_0$, $\tau \in \mathbb{N}^4$ and let $\tau'$ be a permutation of $\tau.$ 
If $t = \begin{bmatrix} t_0 & \tau\end{bmatrix}$ and $t' = \begin{bmatrix} t_0 & \tau'\end{bmatrix}$, then $TE(K^-,t)$ and $TE(K^-, t')$ are cospectral.
\end{lemma}
\begin{proof}
Immediate from Lemma \ref{lemma: charpoly TE Psi}.
\end{proof}
\begin{lemma}\label{lemma: new permute beta part 2}
Let $t_0\in\mathbb{N}_0$ and let  $\tau \in \mathbb{N}^4$ be such that fewer than three entries of $\tau$ are equal. 
Then there exists a $\tau' \not= \tau$, obtained as a permutation of $\tau,$ such that if $t = \begin{bmatrix} t_0 & \tau\end{bmatrix}$ and $t' = \begin{bmatrix} t_0 & \tau'\end{bmatrix}$, then $TE(K^-,t) \not\cong TE(K^-, t')$.
\end{lemma}
\begin{proof}
By contradiction. Recall that $\{1,2\}$ and $\{3,4\}$ are digons in $K^-$, and that two digraphs may only be isomorphic if they contain equal numbers of digons. Let $t = \begin{bmatrix} t_0 & t_1 & t_2 & t_3 & t_4 \end{bmatrix}$, and assume that no three $t_j$'s are equal. ($j=1,\ldots,4$.) Let $t' = \begin{bmatrix} t_0 & t_1 & t_3 & t_2 & t_4 \end{bmatrix}$ and $t'' = \begin{bmatrix} t_0 & t_1 & t_4 & t_3 & t_2 \end{bmatrix}$. Suppose that all three $t$ construct isomorphic expansions. Then by counting digons
\[
    \begin{cases}
        TE(K^-,t) \cong TE(K^-,t') \\
        TE(K^-,t) \cong TE(K^-,t'')
    \end{cases}
    \implies 
    \begin{cases}
    t_1 t_2 + t_3t_4 = t_1t_3 + t_2t_4\\
    t_1 t_2 + t_3t_4 = t_1t_4 + t_2t_3
    \end{cases}
    \implies 
    \begin{cases}
    (t_1 - t_4)(t_2-t_3) = 0 \\
    (t_1-t_3)(t_2-t_4) = 0
    \end{cases}
\]
and thus $t_1 = t_2 = t_3$, $t_1 = t_2 = t_4$, $t_1 = t_3 = t_4$, or $t_2 = t_3 = t_4$, which is a contradiction. 
\end{proof}

From Lemma \ref{lemma: new permute beta part 2}, we find the following necessary condition for SHDS expansions of $K^-$. 
\begin{corollary}\label{cor: neccessity final theorem}
$D = TE\left(K^-,\begin{bmatrix} t_0 & \tau\end{bmatrix}\right)$, $t_0\in\mathbb{N}_0$, $\tau\in\mathbb{N}^4$, is SHDS only if at least three entries of $\tau$ are equal.
\end{corollary}
Intuitively, one might also say that if the condition in Corollary \ref{cor: neccessity final theorem} is satisfied, then any permutation of $t$ would yield digraphs that are isomorphic to one another. In other words, that the condition above is not just necessary, but also sufficient. In the theorem below, which the authors consider to be the main contribution of this paper, we show exactly that. 
\begin{theorem}\label{thm: main theorem}
Let $t$ be an expansion vector. Then $D = TE(K^-,t)$ is SHDS if and only if $t=\begin{bmatrix}t_0 & \tau \end{bmatrix}$, where $t_0\in\mathbb{N}_0$ and $\tau$ is a permutation of $\begin{bmatrix}t_1 & t_1 & t_1 & t_2 \end{bmatrix}$, $t_1,t_2\in\mathbb{N}$. 
\end{theorem}
\begin{proof}Necessity was addressed in Corollary \ref{cor: neccessity final theorem}, so we will only show sufficiency. 
Let $D^*$ be a digraph and suppose that $\Sigma_D = \Sigma_{D^*}$.
If we can show $D \cong D^*$, the claim is true. 
Since $D$ has rank 4 with three positive eigenvalues, we know by Theorem \ref{thm: TR of a digraph with 1 neg ev} that $D^* = TE(K^-,b)$ for some expansion vector $b = \begin{bmatrix} b_0 & b_1 & b_2 & b_3 & b_4\end{bmatrix}$. 
Moreover, by Corollary \ref{cor: explicit spectra}, $D^*$ has eigenvalue $t_1$ with multiplicity 2.
Since $t_1 > 0$, the eigenvectors corresponding to eigenvalue $t_1$ are orthogonal to the nullspace of $H=H(D^*)$. 
Therefore, the eigenvalues are constant on the four nontrivial parts of $D^*$ and 0 on the isolated vertices in $D^*$. 
That is, the eigenvectors $z$ corresponding to $t_1$ satisfy
\begin{equation}z = \begin{bmatrix}\bf{0} & z_1 \bf{j}_1 & z_2 \bf{j}_2 & z_3 \bf{j}_3 & z_4 \bf{j}_4 \end{bmatrix}^\top ~\text{ for }~ z_1, z_2, z_3, z_4 \in \mathbb{C}.\label{eq: v should satisfy}\end{equation}
Recall that Hermitian matrices are diagonalizable \cite{brouwer2011spectra} and that diagonalizable matrices have geometric multiplicities equal to their algebraic multiplicities. Hence, there are two independent eigenvectors $z$ and $w$ that correspond to $t_1$ and satisfy \eqref{eq: v should satisfy}. Moreover, any linear combination of $z$ and $w$ is again an eigenvector for $t_1$. Hence, there is an eigenvector $x$ that is obtained as such a linear combination, which is zero for any one of the nontrivial parts of $D^*$. Suppose that at least one of $x_2,x_3,x_4$ is nonzero and
\begin{equation}x = \begin{bmatrix}\bf{0} & \bf{0} & x_2 \bf{j}_2 & x_3 \bf{j}_3 & x_4 \bf{j}_4 \end{bmatrix}^\top ~\text{ and }~ Hx = t_1 x,\label{eq: x should satisfy}\end{equation}
 where 
\[H = \begin{bmatrix} 0 & 0 & 0 & 0 & 0 \\
                      0 & 0 & J & iJ & -iJ \\
                      0 & J & 0 & -iJ & iJ \\
                      0 & -iJ & iJ & 0 & J \\
                      0 & iJ & -iJ & J & 0 \end{bmatrix} \]
for blocks $H_{ij}$ of appropriate dimensions. Working out the latter equality in \eqref{eq: x should satisfy}, we obtain\\

\typesystem{\systemA}{
              {b_2}x_2                    +               {ib_3}x_3                   -               {ib_4}x_4  =         0,
\textcolor{white}{4x_2} \textcolor{white}{+}             {-ib_3}x_3                   +               {ib_4}x_4  =  {t_1}x_2,
             {ib_2}x_2  \textcolor{white}{+} \textcolor{white}{0x_3}                  +                     x_4  =  {t_1}x_3,
            {-ib_2}x_2                    +                {b_3}x_3 \textcolor{white}{+} \textcolor{white}{0x_4} =  {t_1}x_4}
\systemA\\

\noindent from which it follows that $b_2x_2 = t_1x_2$, $b_3x_3 = t_1x_3$ and $b_4x_4 = t_1x_4$. Moreover, at most one of $x_2,x_3,x_4$ can be zero and therefore at least two of $b_2,b_3,b_4$ are equal to $t_1$. But similarly, this implies that at least two among each three of $b_1,b_2,b_3,b_4$ equals $t_1$, and thus at least three of $b_1,b_2,b_3,b_4$ are equal to $t_1.$ Thus, since the last four numbers may be permuted to any order without changing the spectrum, we may assume without loss of generality that $b = \begin{bmatrix} b_0& t_1 & t_1 & t_1 & b_4 \end{bmatrix}$. 


Finally, by comparing the number of edges in the underlying graph, we obtain $b_4 = t_2$, and since $|\Sigma_D| = |\Sigma_{D^*}|$, we have $b_0 = t_0$. Hence, if $\Sigma_D = \Sigma_{D^*}$ then $b = t$ and thus $D^* \cong D$.
\end{proof}
Thus, by Theorem \ref{thm: main theorem}, we may be certain that we may uniquely determine each digraph whose Hermitian spectrum is of the form 
\[\left\{-\sqrt{3t_1t_2 + t_1^2}-t_1,~\sqrt{3t_1t_2+t_1^2}-t_1,~t_1^{[2]},~0^{[t_0 + 3t_1 + t_2 - 4]} \right\}, t_0\in\mathbb{N}_0, t_1,t_2\in\mathbb{N}.\]

\section{Closing remarks regarding WHDS digraphs}
While most of this paper has been concerned with strong determination by the Hermitian spectrum, we conclude with some remarks on its weaker counterpart. It stands to reason that the digraphs that are not SHDS due to the problem illustrated in Example \ref{ex2} might be WHDS, as most of the machinery that was used to show Theorem \ref{thm: main theorem} is still in place. We first observe that it is in general not true that any expansion of $K^-$ is WHDS, after which we prove an analogue to Proposition \ref{prop: wang characterization}. 

As mentioned before, relaxing the connectivity assumption may yield cospectral digraphs with distinct underlying graphs, which implies that said digraphs are not switching equivalent. Recall the following example by \citet{wang2017mixed}.
\example{ex wang}{ \cite{wang2017mixed} Let $D = TE\left(T^-, \begin{bmatrix} t_0 & 3 & 3 & 18\end{bmatrix}\right)$  and $D^* = TE\left(T^-, \begin{bmatrix} t_0+4 & 2 & 9 & 9\end{bmatrix}\right)$, for some $t_0\in\mathbb{N}_0$. Then $\Sigma_D = \Sigma_{D^*}$ (plug in \eqref{eq: charpoly TE Psi}) regardless of $t_0$. Since $D$ and $D^*$ do not contain an equal number of isolated vertices, they cannot be switching equivalent. Thus, there are infinitely many pairs of cospectral mates. Most notably, suppose that $t_0=0$. Then $D$ is connected, but not switching equivalent to $D^*$, while they are cospectral. Thus, $D$ is not WHDS, even when $t_0 = 0$.}

As we have seen throughout this paper, there are many parallels between $T^-$ and $K^-$. This has caused us to believe that a similar phenomemon occurs for the negative tetrahedron. 
As we will see shortly, there are indeed pairs of expansion vectors $t,t'$ for $K^-$ such that $D=TE(K^-,t)$ and $D'=TE(K^-,t')$ are both connected and have the same nonzero eigenvalues, but not the same number of vertices, which allows us to formulate an analogue of Proposition \ref{prop: wang characterizaton}. 
\begin{proposition}\label{prop: inf many cospectral pairs}
There are infinitely many rank $4$ digraphs with exactly one negative eigenvalue that are not WHDS.
\end{proposition}
The correctness of Proposition \ref{prop: inf many cospectral pairs} follows immediately from the following example. 
\example{ex: TE K-}{ Let $t_0\in\mathbb{N}_0$ and let $t = \begin{bmatrix} t_0 & 9 & 18 & 20 & 60\end{bmatrix}$ and  $t' = \begin{bmatrix} t_0+4 & 10 & 12 & 36 & 45\end{bmatrix}$. Then $D = TE(K^-,t)$ and $D' = TE(K^-,t')$ both have characteristic polynomial 
\begin{equation}\label{eq: charpoly example}\chi(\mu) = \mu^{103+t_0}(\mu^4 - 3522 \mu^2 + 90720 \mu - 583200).\end{equation}
Moreover, as the two contain distinct numbers of isolated vertices, they are clearly not switching equivalent.}

However, it should also be noted that such examples have proven to be extremely rare. 
If one enumerates all vectors $t=\begin{bmatrix}0 & t_1 & t_2 & t_3 & t_4\end{bmatrix}$ with\footnote{Note that this does not restrict the investigation, as we are not interested in switching equivalent pairs.} $0<t_1\leq t_2 \leq t_3 \leq t_4 \leq 104$ and computes the characteristic polynomials corresponding to $TE(K^-,t)$ via Lemma \ref{lemma: charpoly TE Psi}, we find as few as five characteristic polynomials whose nonzero roots do not occur uniquely. 
Out of these five polynomials, $\chi(\mu)$ as in \eqref{eq: charpoly example} is of smallest order; the corresponding digraphs contain at least 107 vertices.
For the remaining digraphs, we obtain strong evidence that they are WHDS.
In particular, we may thus conclude the following.
\begin{proposition}
Any digraph of order less than $107$ that has rank $4$ and exactly one negative eigenvalue is WHDS.
\end{proposition}
\begin{proof}
It is clear that any pair of digraphs $D = TE(K^-,[0~~\tau])$, $D'= TE(K^-,[0~~\tau'])$ with $|V(D)|, |V(D')| < 107$ and equal nonzero eigenvalues would have occurred in the performed enumeration. Since they did not, and the number of isolated vertices that may be added is bounded by the assumption on the order of the digraphs, the claim follows. 
\end{proof}
\begin{proposition}
Let $t_0\in\mathbb{N}_0, \tau\in\mathbb{N}$ with $\tau_j\leq 7 $, $j=1,\ldots,4$. Then $TE(K^-,[t_0~~\tau])$ is WHDS. 
\end{proposition}
\begin{proof} Let $t = [t_0~~\tau]$ with $\tau_j\leq 7$, $j\in[4]$, and $t' = [t_0~~\tau']$ where $\tau'$ contains at least one element larger than 100. Then
\[\sum_{1\leq i < j \leq 4} t_it_j \leq 294 < 306 \leq \sum_{1\leq i < j \leq 4} t'_it'_j.\]
Thus, if there were expansion vectors that would have yielded the same characteristic polynomial coefficients as $t$,  they would have occurred in the performed enumeration. 
\end{proof}

Moreover, if we assume connectivity, we obtain a result similar to Proposition \ref{prop: wang characterization}.
\begin{proposition}\label{prop: connected cospectral rank 4 digraphs are sw eq}
Any two connected, cospectral rank $4$ digraphs with exactly one negative eigenvalue are switching equivalent.
\end{proposition}
\begin{proof}
Let $D_1$ and $D_2$ be connected, cospectral, rank 4 digraphs with exactly one negative eigenvalue. By Theorem \ref{thm: TR of a digraph with 1 neg ev}, $TR(D_1) = TR(D_2) = K^-$. Let $t,s$ be such that $D_1 = TE(K^-,t)$ and $D_2 = TE(K^-,s)$. (Note that $t_0 = s_0 = 0$ is implied.) 

We first observe that $s$ may be obtained as a permutation of $t$. Indeed, consider equations \eqref{eq: whds proof 1}-\eqref{eq: whds proof 4}, which should all hold simultaneously by Lemma \ref{lemma: charpoly TE Psi}. (Note that \eqref{eq: whds proof 1} is implied by $t_0 = s_0$.) Then, if these elementary symmetric polynomials in $\{t_1,t_2,t_3,t_4\}$ and $\{s_1,s_2,s_3,s_4\}$ are the same, $s$ is indeed distinct from $t$ by at most a reordering of its entries. 
\startCustomModel
\begin{align}
    \sum_{j=1}^4 t_j &= \sum_{j=1}^4 s_j \label{eq: whds proof 1}\\
    \sum_{1\leq i < j \leq 4} t_i t_j &= \sum_{1\leq i < j \leq 4}  s_is_j\\
    \sum_{1\leq i < j < k\leq 4} t_it_jt_k &= \sum_{1\leq i < j < k\leq 4} s_is_js_k\\
    \prod_{j=1}^4 t_j &= \prod_{j=1}^4 s_j \label{eq: whds proof 4}
\end{align}
\resetEquation

Now, observe that any permutation of $t$ might be obtained by a sequence of pairwise interchanges of elements. 
Thus, if we can show that every such pairwise exchange yields a switching equivalent digraph, we are done. 
Suppose that $t'$ is obtained from $t$ by interchanging elements $t_u$ and $t_v$, $u,v\in [4],$ and let $D = TE(K^-,t)$ and $D'= TE(K^-,t')$.
We distinguish two cases: either $\{u,v\}$ is a digon in $K^-$ or it is not. 
Note that we may assume without loss of generality that $(u,v) = (1,2)$ in the former case, and $(u,v) = (1,3)$ in the latter.

Suppose that $(u,v) = (1,2)$. Then $D'\cong D^\top$, and thus $D$ is switching equivalent to $D'$. Finally, suppose that $(u,v) = (1,3)$ and let $S = \text{Diag}\left(\begin{bmatrix}-i\cdot \mathbf{j}_{t_1} & \mathbf{j}_{t_2} & i\cdot \mathbf{j}_{t_3} & \mathbf{j}_{t_4} \end{bmatrix}\right)$. 
If we let $D_S$ be the digraph determined by the Hermitian $H_S = S^{-1}H(D)S$, then $D_S \cong D'$. Since $D_S$ was obtained from $D$ by a four-way switching, $D$ is again switching equivalent to $D'$, which completes the proof. 
\end{proof}

\scriptsize
\bibliography{mybib}{}
\bibliographystyle{abbrvnat}
\normalsize


\appendix

\section{Fraction of self-converse digraphs}
\label{sec: fraction self-converse}
In the interest of completeness, we include the following numerical evidence, that illustrates the scarceness of SHDS digraphs, which necessarily must be self-converse. 

\begin{conjecture}\label{conj: fraction self converse to zero} Let $\mathcal{SC}_n$ be the number of self-converse digraphs on $n$ vertices, and let $\mathcal{D}_n$ be the number of digraphs on $n$ vertices. Then:
\begin{equation}
\frac{\mathcal{SC}_n}{\mathcal{D}_n} \to 0 \text{ as } n\to\infty.
\end{equation}
\end{conjecture}

\noindent \textit{Numerical evidence for Conjecture \ref{conj: fraction self converse to zero}.} ~Define $f(n)$ to be $f(n) = d^{sc}_n(1)/d_n(1)$, where $d^{sc}_n$ and $d_n$ are the counting polynomials that count respectively the number of self-converse digraphs on $n$ vertices and the number of digraphs on $n$ vertices, as defined in \cite{harary1966enumeration, harary1955number}. 
Conjecture \ref{conj: fraction self converse to zero} claims that $f(n)\to 0$ as $n\to \infty$. 
Since explicit formulas are available, one may numerically evaluate them; the results are shown in Table \ref{tab: fraction self converse}. 
By the shear speed at which the fraction of self-converse graphs converges to zero, one is easily convinced that the claim is true. 

\begin{table}[h!]
\begin{center}
\begin{tabular}{ccccccccccc}
$n$  & 3 & 4 & 5& 6& 7 &8 \\
$f(n)$ & 6.25$\cdot 10^{-1}$&   3.21$\cdot 10^{-1}$&   7.36$\cdot 10^{-2}$&   9.87$\cdot 10^{-3}$&   6.16$\cdot 10^{-4}$&   2.20$\cdot 10^{-5}$\\  
\hline
$n$ &9 &10  & 11 & 12 & 13 & 14  \\
$f(n)$ &   3.89$\cdot 10^{-7}$&   3.79$\cdot 10^{-9}$ & 1.85$\cdot 10^{-11}$&   4.89$\cdot 10^{-14}$&   6.50$\cdot 10^{-17}$&   4.58$\cdot 10^{-20}$\\
\hline
$n$  & 15& 16& 17 &18 &19 &20\\
$f(n)$ &   1.63$\cdot 10^{-23}$&   3.06$\cdot 10^{-27}$&   2.90$\cdot 10^{-31}$&   1.43$\cdot 10^{-35}$&   3.59$\cdot 10^{-40}$&   4.64$\cdot 10^{-45}$
   \end{tabular}
\caption{\label{tab: fraction self converse} The fraction $f(n)$ of digraphs of size $n$ that is self-converse.}
\end{center}
\end{table}


\end{document}

%% file: commands.tex
\newcounter{example}
\newenvironment{example}[2]{
\refstepcounter{example}
\label{#1}\begin{leftbar}\textbf{\noindent \hspace{-5.3mm}Example \theexample}.#2\end{leftbar}
}{}

\newcounter{drafts}

\newcounter{module}
\newcounter{model}
\newcommand{\startCustomModel}{
\newcounter{tmp}
\setcounter{tmp}{\value{equation}}
\refstepcounter{tmp}
\setcounter{equation}{0}
\renewcommand{\theequation}{\thetmp\alph{equation}}
}

\newcommand{\resetEquation}{
\renewcommand{\theequation}{\arabic{equation}}
\setcounter{equation}{\value{tmp}}
}

\newcommand \linedabstractkw[2]{
  \renewcommand\maketitlehookd{%
    \mbox{}\medskip\par
    \centering
    \hrule\medskip
    \begin{minipage}{0.9\textwidth}
    \textbf{Abstract}\\ #1\\
    
    \textit{Keywords: }#2
    \end{minipage}\medskip\hrule\medskip
    }      
}

\newcommand \ERPWtemplate[3]{
\usepackage{fancyhdr}
\pagestyle{fancy}
\fancyhf{}
\fancyhead[RO]{#1}
\fancyhead[LO]{#2}
\fancyfoot[CO]{\thepage}

\renewcommand{\headrulewidth}{1pt}
\renewcommand{\footrulewidth}{1pt}
}